\documentclass[12pt,a4paper]{article}
\usepackage[top=20mm,bottom=20mm]{geometry}
\usepackage{amsmath,defns,euler,natbib,url,graphicx,microtype}
\usepackage{amsthm,amsfonts}
\usepackage{footmisc}

\usepackage{hyperref} 
\hypersetup{pdftex,
            backref=true,
            hyperindex=true,
            colorlinks=true,
            citecolor=blue,
            bookmarks=true,
            breaklinks=true}

\newcommand{\uu}{{\bar u}}

\renewcommand{\vec}[1]{\text{\boldmath$#1$}}

\newcommand{\shift}{\text{\small$\mathcal E$}}

\ifx\de\undefined
\newcommand{\de}[2]{\mathchoice
  {\frac{d #2}{d #1}}% display
  {{d #2}/{d #1}}% text
  {{d #2}/{d #1}}% script
  {{d #2}/{d #1}}% scriptscript
  }
\fi
\newcounter{i}
\newcommand{\HH}{\ensuremath{\mathbb H}}

\newcommand{\MM}{\ensuremath{\mathbb M}}
\newcommand{\cC}{\ensuremath{\mathcal C}}
\newcommand{\cM}{\ensuremath{\mathcal M}}
\newcommand{\cL}{\ensuremath{\mathcal L}}
\newcommand{\cK}{\ensuremath{\mathcal K}}

\allowdisplaybreaks

\title{Multiscale modelling couples patches of nonlinear wave-like simulations}

\author{Meng Cao\thanks{School of Mathematical Sciences,
University of Adelaide, South Australia 5005.  \protect\url{mailto:meng.cao@adelaide.edu.au} or \protect\url{mailto:mengcao1188216@gmail.com}}
\and 
A.~J. Roberts\thanks{School of Mathematical Sciences,
University of Adelaide, South Australia 5005.  \protect\url{mailto:anthony.roberts@adelaide.edu.au}}}

\date{\today}

\begin{document}
    
\maketitle

\begin{abstract}
The multiscale gap-tooth scheme is built from given microscale simulations of complicated physical processes to empower macroscale simulations.
By coupling small patches of simulations over unsimulated physical gaps, large savings in computational time are possible. 
So far the gap-tooth scheme has been developed for dissipative systems, but wave systems are also of great interest.
This article develops the gap-tooth scheme to the case of nonlinear microscale simulations of wave-like systems. 
Classic macroscale interpolation provides a generic coupling between patches that achieves arbitrarily high order consistency between the multiscale scheme and the underlying microscale dynamics.
Eigen-analysis indicates that the resultant gap-tooth scheme empowers feasible computation of large scale simulations of wave-like dynamics with complicated underlying physics.
As an pilot study, we implement numerical simulations of dam-breaking waves by the gap-tooth scheme.
Comparison between a gap-tooth simulation, a microscale simulation over the whole domain, and some published experimental data on dam breaking, demonstrates that the gap-tooth scheme feasibly computes large scale wave-like dynamics with computational savings.
\end{abstract}

\tableofcontents

\section{Introduction}
\label{sec:intro}

Mathematical equations describing floods and tsunamis are typically written at the macroscale of kilometres. 
But the level at which the underlying turbulent fluid physics are best understood is at the much finer sub-metre scale. 
Although macroscale models for floods and tsunamis are well established, for many multiscale and multiphysics wave-like problems good macroscale descriptions (good closures) do not exist. 
We aim to empower scientists and engineers to use brief bursts of microscale wave-like simulation on small patches of the space-time domain in order to make efficient accurate macroscale simulations without ever knowing a macroscale closure.

Many multiscale modelling techniques have been developed for dissipative systems~\cite[e.g.]{E:2003fk, Kevrekidis:2003fk, Roberts:2005uq, Hou:2008fk}.
Our macroscopic modelling further develops the equation-free gap-tooth scheme \cite[e.g.]{Gear03, Samaey03b, Samaey:2005fk, Samaey2009} to empower simulation of wave-like systems over large time and space scales from a \emph{given} microscopic simulator.
We suppose that the wave-like microscale simulator is computationally expensive so that only small time and spatial domain simulations are feasible: one example of future interest is direct numerical simulation of depth resolved turbulent fluid floods.
The microscale simulator provides the necessary data for the macroscopic model, so whenever the microscale simulator improves, then the overall macroscale simulation will correspondingly improve.

\begin{figure}
\begin{center}
\setlength{\unitlength}{0.7ex}
\begin{picture}(120,70)(0,0)
\put(15,0){\includegraphics{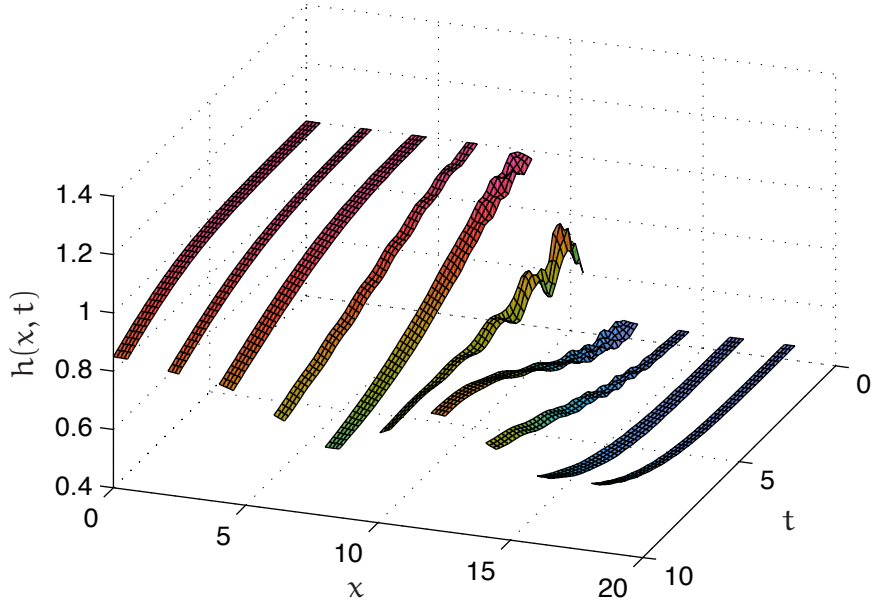}}
\put(10,0){\rotatebox{90}{\hspace{16ex}$h(x,t)$}}
\put(45,0){$x$}
\put(90,7){$t$}
\end{picture}
\end{center}
\caption{Indicative gap-tooth simulation of a dam-break shows the water depth~\(h(x,t)\): microscale computations are only done in small patches of the spatial domain~\(x\), but here over all time; at time zero (back) the dam breaks and at later times (front) the water forms a turbulent bore that propagates to the right.}
\label{fig:PatchDamInt}
\end{figure}

This article develops and theoretically supports the gap-tooth method for general wave-like systems, and as an indicative application and test, applies the methodology to a Smagorinski model of turbulent shallow water flow.
For example, Figure~\ref{fig:PatchDamInt} exhibits the indicative gap-tooth simulation of a dam-break showing the water depth~\(h(x,t)\), where  the microscale computations are only done in small patches of the spatial domain~\(x\) over all time.
Here the scheme uses microscale simulations of shallow water flow on small patches of space (Section~\ref{sec:micro}), coupling the simulations over the intervening space, to simulate floods over a macroscale.
Potential future applications could improve modelling of sediment erosion, transport, and deposition.
The gap-tooth method and our theoretical support (Section~\ref{sec:couple}) adapts to whatever microscale simulator is provided. 
The analysis of Section~\ref{sec:couple} indicates that the patches can occupy as small a fraction of space as is necessary for a good microscale simulation without affecting macroscale accuracy, thus indicating large computational gains are feasible with the methodology.
%Our early work~\cite[e.g.]{Cao2013} simulated the linear wave-like system using this gap-tooth scheme.
%\ajr{Clarify a little what was achieved.  e.g. high-order consistency to basic wave equation.}
%Numerical simulations and eigenvalue analysis showed this gap-tooth scheme works well with the linear microscale simulators in the wave-like system.

Previous research \cite[]{Cao2013} explored classic linear non-dispersive waves and found classic interpolation on a macroscale staggered grid ensured ensure high order consistency to the linear wave equation.
Section~\ref{sec:couple} proves for general linear microscale wave systems that the patch coupling condition~\eqref{patch:cph} ensures arbitrarily high order consistency between the gap-tooth scheme and the underlying microscale dynamics, and also establishes consistency for a class of nonlinear wave systems.
Further, section~\ref{sec:nsm} discusses nonlinear wave-like systems further and establishes the gap-tooth scheme with patches coupled with~\eqref{patch:cph} have a slow manifold \cite[e.g.]{Boyd95, MacKay03} that forms the macroscale dynamics.
This article focusses on waves in one spatial dimension with the expectation that generalisation to multiple space dimensions will be analogous to that for dissipative systems \cite[e.g.]{Roberts2011a}.

%Section~\ref{sec:micro} establishes a gap-tooth multiscale simulation of nonlinear wave dynamics on a staggered grid.
%Then section~\ref{sec:couple} introduces the developed coupling conditions by interpolating the macroscale grid values to provide boundary values on each patch.
%These coupling conditions achieve arbitrarily high order consistency between the patch scheme and the underlying microscale dynamics for general microscale systems.

Section~\ref{sec:Eig} implements the gap-tooth method by coupling small patches of the given microscale simulations of the shallow water Smagorinski model described in Section~\ref{sec:micro}.
Numerical eigenvalue analysis supports the theoretical results that there is an appropriate slow manifold of the macroscale dynamics in this application of the gap-tooth scheme.

Numerical simulations show that the gap-tooth coupling condition~\eqref{patch:cph} works well for a range of shallow water flows.
Section~\ref{sec:dam} applies the gap-tooth simulation to a dam-break (Figure~\ref{fig:PatchDamInt}), then compares the gap-tooth simulation with the microscale simulation over the whole domain, and with some experimental data of \cite{Stansby1998}.

\section{The nonlinear microscale water wave model}
\label{sec:micro}

This section describes the nonlinear microscale simulator of the nonlinear shallow water wave \pde\ derived from the Smagorinski model of turbulent flow \cite[]{Roberts:2008fk, Cao2014}.
Often, wave-like systems are written in terms of two conjugate variables, for example, position and momentum density, electric and magnetic fields, and water depth~\(h(x,t)\) and mean lateral velocity~\(\uu(x,t)\) as herein.
This article uses the example of shallow water waves, but applies to any wave-like system in the form
\begin{equation}
\D th=-c_1\D x\uu+f_1[h,\uu]
\quad\text{and}\quad
\D t\uu=-c_2\D xh+f_2[h,\uu],
\label{eq:genwaveqn}
\end{equation}
where the brackets indicate that the nonlinear functions~\(f_\ell\) may involve various spatial derivatives of the fields~\(h(x,t)\) and ~\(\uu(x,t)\).
Specifically, this section invokes a nonlinear Smagorinski model of turbulent shallow water \cite[e.g.]{Roberts:2008fk, Cao:2012fk} along an inclined flat bed: let $x$~measure position along the bed and in terms of fluid depth~$h(x,t)$ and depth-averaged lateral velocity~$\uu(x,t)$ the model \pde{}s are
\begin{subequations}\label{eqs:patch:N}%
\begin{align}
\frac{\partial h}{\partial t}&=-\frac{\partial(h\uu)}{\partial x}\,,\label{patch:Nh}
\\
\frac{\partial \uu}{\partial t}&={}0.985\left(\tan\theta-\frac{\partial h}{\partial x}\right)-0.003\frac{\uu|\uu|}{h}-1.045\uu\frac{\partial \uu}{\partial x}+0.26h|\uu|\frac{\partial^2\uu}{\partial x^2}\,,\label{patch:Nu}
\end{align}
\end{subequations}
where~$\tan\theta$ is the slope of the bed.
Equation~\eqref{patch:Nh} represents conservation of the fluid.
The momentum \pde~\eqref{patch:Nu} represents  the effects of turbulent bed drag~$\uu|\uu|/h$, self-advection~$\uu\D x\uu$, nonlinear turbulent dispersion~$h|\uu|\DD x\uu$, and gravitational hydrostatic forcing~$\tan\theta-\D xh$.

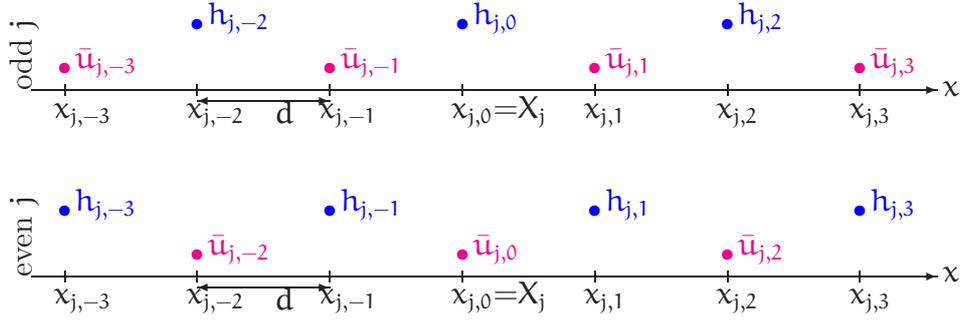
\begin{figure}
\begin{center}
\setlength{\unitlength}{0.8ex}
\begin{picture}(86,40)
%\put(0,0){\framebox(86,28){}}
\put(0,17){%shift upper part upwards
\put(0,3){\rotatebox{90}{odd \(j\)}}
\put(2,3){\vector(1,0){82}}
\put(84.5,2.5){$x$}
%\newcounter{i}
\setcounter{i}{-4}
\multiput(5,2.5)(12,0){7}{\line(0,1){1}%
  \stepcounter{i}%
  \put(-1,-2){$x_{j,\arabic{i}}
  \ifnum\arabic{i}=0{=}X_j\fi$}
  }
\put(23,2){\vector(-1,0){6}
  \put(6,0){\vector(1,0){6}}
  \put(7,-2){$d$}
  }
\setcounter{i}{-4}
\multiput(5,7)(12,0){7}{%
  \stepcounter{i}%
  \ifodd\arabic{i}\color{magenta}%
    \put(0,-2){%
    \circle*{1}$\uu_{j,\arabic{i}}$
    }
  \else\color{blue}%
    \put(0,2){%
    \circle*{1}$h_{j,\arabic{i}}$
    }
  \fi
}
}% end of upper part
\put(0,3){\rotatebox{90}{even \(j\)}}
\put(2,3){\vector(1,0){82}}
\put(84.5,2.5){$x$}
%\newcounter{i}
\setcounter{i}{-4}
\multiput(5,2.5)(12,0){7}{\line(0,1){1}%
  \stepcounter{i}%
  \put(-1,-2){$x_{j,\arabic{i}}
  \ifnum\arabic{i}=0{=}X_j\fi$}
  }
\put(23,2){\vector(-1,0){6}
  \put(6,0){\vector(1,0){6}}
  \put(7,-2){$d$}
  }
\setcounter{i}{-4}
\multiput(5,7)(12,0){7}{%
  \stepcounter{i}%
  \ifodd\arabic{i}\color{blue}%
    \put(0,2){%
    \circle*{1}$h_{j,\arabic{i}}$
    }
  \else\color{magenta}%
    \put(0,-2){%
    \circle*{1}$\uu_{j,\arabic{i}}$
    }
  \fi
}
\end{picture}
\end{center}
\caption{Scheme of the staggered grid points of the depth $h_{j,i}$~(blue points) and velocity $\uu_{j,i}$~(magenta points) at the \(i\)th~micro-grid point on the odd $j$th patch~(top) and the even $j$th patch~(bottom).  This diagram shows the cases for \(n=5\) interior grid points in each patch.}
\label{patchMicro}
\end{figure}

In practice, the microscale simulator will typically be either a spatial discretisation such as finite difference \cite[e.g.]{Bijvelds99}, finite element, or finite volume \cite[e.g.]{LeVeque2011}, or a particle based method such as lattice Boltzmann \cite[e.g.]{Liu2009b}, molecular dynamics \cite[e.g.]{Southern2008}, or smoothed-particle hydrodynamics \cite[e.g.]{Monaghan92}.
%For definiteness, we use finite differences to approximate the \textsc{pde}s~\eqref{eq:genwaveqn} by a set of \textsc{ode}s on a microscale staggered grid.
Our microscale simulator for the nonlinear \pde{}s~\eqref{eqs:patch:N} is a spatial discretisation on a fine-scale staggered grid within each patch.
Figure~\ref{patchMicro} shows the staggered grid points of the depth~$h$ and depth-averaged velocity~$\uu$ on a patch.
Because we propose that the macroscale gap-tooth scheme employ a macroscale staggered grid, there are two types of alternating microscale patches corresponding to even and odd macroscale index~\(j\).
In the~$j$th patch, define a microscale staggered grid of spacing~$d$: at the~$i$th point of the micro-grid of the $j$th~patch define  the depth~$h_{j,i}$ when \(j-i\)~is odd, and define the depth-averaged lateral velocity~$\uu_{j,i}$ when \(j-i\)~is even (Figure~\ref{patchMicro}).
Thus, approximate the Smagorinski shallow water \pde{}s~\eqref{eqs:patch:N} on the~$j$th patch with centred differences in microscale space as the discrete
\begin{subequations}\label{eqs:patch:Nji}%
\begin{align}
\frac{\partial h_{j,i}}{\partial t}&={}-\frac{(h_{j,i+2}+h_{j,i})\uu_{j,i+1}}{4d}+\frac{(h_{j,i-2}+h_{j,i})\uu_{j,i-1}}{4d}\,,\label{patch:Nhji}
\\
\frac{\partial \uu_{j,i}}{\partial t}&={}0.985\left(\tan\theta-\frac{h_{j,i+1}-h_{j,i-1}}{2d}\right)-0.003\frac{\uu_{j,i}|\uu_{j,i}|}{h_{j,i}}
\nonumber\\&\quad
{}-1.045\uu_{j,i}\frac{\uu_{j,i+2}-\uu_{j,i-2}}{4d}
+0.26h_{j,i}|\uu_{j,i}|\frac{\uu_{j,i+2}-2\uu_{j,i}+\uu_{j,i-2}}{4d^2}\,,\label{patch:Nuji}
\end{align}
\end{subequations}
where~$d$ is the microscale spatial step within a patch.
Such a microscale lattice simulator is known to be consistent with the \pde{}s~\eqref{eqs:patch:N} to a error~\Ord{d^2} which typically is negligible for small patches.

The dam break shown in Figure~\ref{fig:PatchDamInt} was generated by such microscale simulations in patches coupled together across unresolved space. 
Time integration was done by Matlab \verb|ode15s|.

\begin{figure}
\begin{center}
\setlength{\unitlength}{1ex}
\begin{picture}(74,20)
%\put(0,0){\framebox(74,15){}}
\put(8,0){% apply shift right to all
\put(0,3){\vector(1,0){60}}
\put(60.5,2.5){$x$}
\setcounter{i}{0}
\multiput(5,2.5)(12,0){5}{\line(0,1){1}%
  \stepcounter{i}%
  \put(-1,-2){$X_{\arabic{i}}$}
  }
\put(23,2){\vector(-1,0){6}
  \put(6,0){\vector(1,0){6}}
  \put(5,-2){$D$}
  }
\setcounter{i}{0}
\multiput(5,9)(12,0){5}{%
  \stepcounter{i}%
  \ifodd\arabic{i}\color{blue}%
    \def\myv{+}\def\myU{H}\def\myu{h}\def\myh{u}%
  \else\color{magenta}%
    \def\myv{-}\def\myU{U}\def\myu{u}\def\myh{h}%
  \fi%
  \put(0,\myv3){\circle*{1}}
  \put(-1,\myv5){\put(0,-0.5){$\myU_{\arabic{i}}$}}
  \put(-2,\myv3){\thicklines%
    \line(1,0)4\,$\myu_{\arabic{i}}(x,t)$
    }
  \put(-2,\myv1.5){\circle{0.6}}
  \put(+2,\myv1.5){\circle{0.6}}
  \put(-2,\myv1.5){\thicklines%
    \line(1,0)4\,$\myh_{\arabic{i}}(x,t)$
    }
  }
\thinlines
\setcounter{i}{0}
\multiput(5,9)(12,0){5}{%
  \stepcounter{i}%
  \ifodd\arabic{i}\color{cyan}%
    \def\myv{+}\def\myU{H}\def\myu{h}\def\myh{u}%
  \else\color{green}%
    \def\myv{-}\def\myU{U}\def\myu{u}\def\myh{h}%
  \fi%
  \put(0.5,\myv3){\vector(2,-\myv1){9}}
  \put(0.5,\myv3){\vector(3,-\myv1){13}}
  \put(-0.5,\myv3){\vector(-2,-\myv1){9}}
  \put(-0.5,\myv3){\vector(-3,-\myv1){13}}
  }
}%end shift right
\end{picture}
\end{center}
\caption{The macroscale scheme interpolates macroscale grid values~$H_j$ and~$U_j$ to provide edge values on each patch.
The green arrows provide edge values of odd~$j$ patches by interpolating macroscale grid values~$U_j$.
The cyan arrows provide edge values of even~$j$ patches by interpolating macroscale grid values~$H_j$.}
\label{patchMacro}
\end{figure}
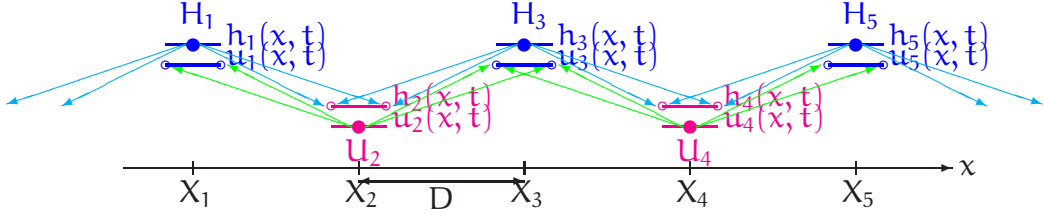

The macroscale grid is also staggered: Figure~\ref{patchMacro} shows the alternating patches.
%\cite{Cao2013} proposed coupling conditions that connects these patches by interpolating macroscale grid values to provide boundary values for each microscale patch simulation, as shown schematically in Figure~\ref{patchMacro}.
Consequently, the macroscale model is to be parametrised by the macroscale grid values
\begin{equation}
U_j(t):=\uu_j(X_j,t)\text{ for even \(j\), and}\quad
H_j(t):=h_j(X_j,t)\text{ for odd }j,
\label{patch:macrogrid}
\end{equation}
as shown in Figure~\ref{patchMacro}.
Following the pilot study of \cite{Cao2013}, for the specific analysis and simulations of Sections~\ref{sec:Eig} and~\ref{sec:dam}  cubic interpolation of macroscale values from the nearest three patches on either side of a patch provided boundary values on that patch.  
That is, on patches with odd index~\(j\) the boundary value of the microscale \(\uu\)-field are determined as
\begin{subequations}\label{eqs:cubic}%
\begin{equation}
\uu_j(X_j\pm rD,t)=
\left[\mu\pm\tfrac12r\delta+\tfrac18(-1+r^2)\mu\delta^2\pm\tfrac{1}{48}(-r+r^3)\delta^3\right]U_j(t),
\label{eq:cubicu}
\end{equation}
for centred mean and difference operators, \(\mu U_j:=(U_{j+1}-U_{j-1})/2\) and \(\delta U_j:=U_{j+1}-U_{j-1}\) respectively (the staggered macroscale grid of Figure~\ref{patchMacro} requires these non-standard definitions).
Correspondingly, on patches with even index~\(j\) the boundary value of the microscale \(h\)-field are 
\begin{equation}
h_j(X_j\pm rD,t)=
\left[\mu\pm\tfrac12r\delta+\tfrac18(-1+r^2)\mu\delta^2\pm\tfrac{1}{48}(-r+r^3)\delta^3\right]H_j(t).
\label{eq:cubich}
\end{equation}
\end{subequations}
The expansion~\eqref{eq:shif} justifies that these are cubic interpolation of the four surrounding macroscale grid values.
This interpolation couples the macroscale staggered grid of patches across un-simulated space, as seen in the dam break of Figure~\ref{fig:PatchDamInt}, to form a well-posed simulation.

\section{Numerical simulations and eigenvalues show separation of scales}
\label{sec:Eig}

We numerically explore the macroscale turbulent fluid flow on a slightly inclined flat bed using the gap-tooth scheme with the nonlinear microscale simulator~\eqref{eqs:patch:Nji} and the cubic coupling conditions~\eqref{eqs:cubic}.
%Section~\ref{sec:couple} shows the high order consistency of the coupling conditions~\eqref{patch:cph}.
%In these numerical simulations, the coupling conditions~\eqref{patch:cph} are truncated to errors~$\mathcal O(\gamma^3)$ to give a cubic interpolation of the neighbouring macroscale values to approximate the edge values of a patch.
Numerical simulations are straightforwardly implemented for equations \eqref{eqs:patch:Nji} on staggered grids in space of Figures~\ref{patchMicro} and~\ref{patchMacro}.
Time integration was performed by Matlab's \verb|ode15s|.

\begin{figure}
\centering
\begin{tabular}{r@{\ }c}
 $t=0$\\
\rotatebox{90}{\hspace{7ex}$h$ and $\uu$ }&
\includegraphics[height=23ex]{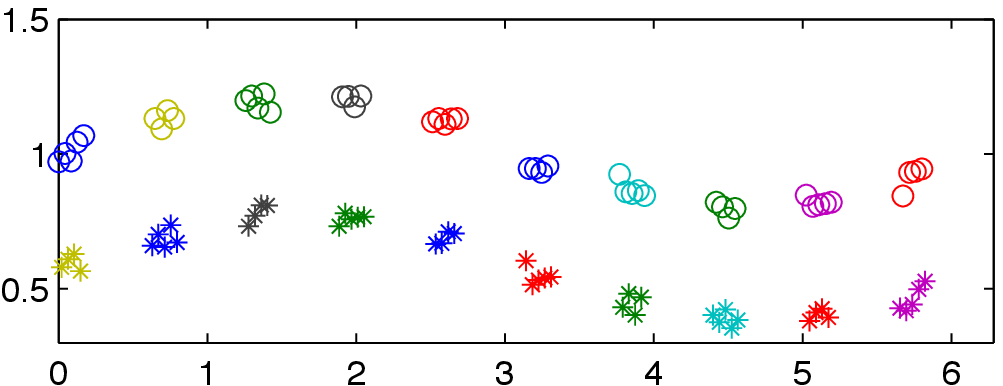}\\
& $x$ \\
$t=2$\\
\rotatebox{90}{\hspace{7ex}$h$ and $\uu$ }&
\includegraphics[height=23ex]{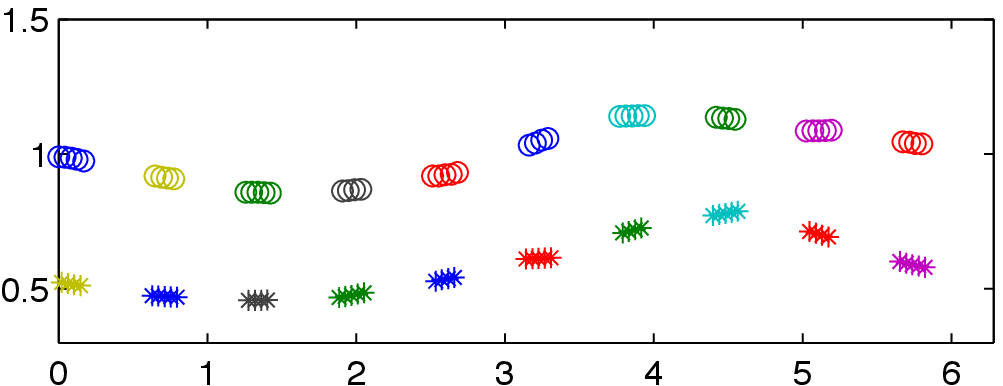}\\
& $x$ \\
$t=4$\\
\rotatebox{90}{\hspace{7ex}$h$ and $\uu$  }&
\includegraphics[height=23ex]{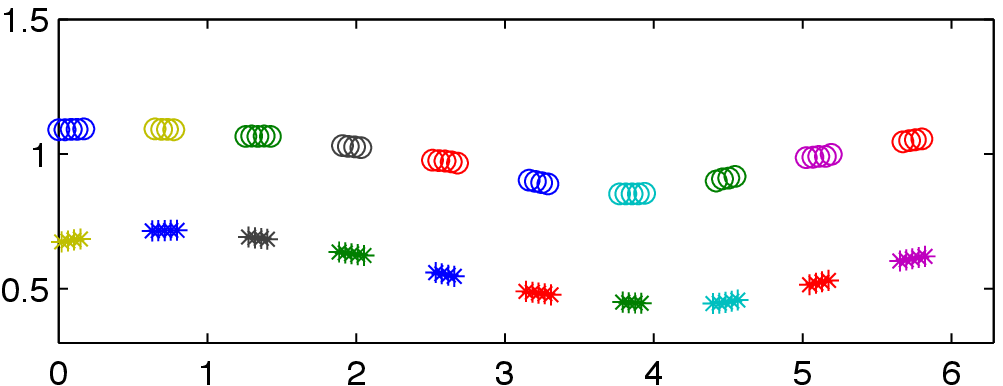}\\
& $x$
\end{tabular}
\caption{Gap-tooth simulation of water depth $h$~(circles) and average lateral velocity $\uu$~(stars) by the nonlinear microscale simulator~\eqref{eqs:patch:Nji} and the coupling conditions~\eqref{eqs:cubic} on doman $\Omega=[0,2\pi]$ via $m=10$ patches and $n=9$ microscale grid points on each patch at three times. 
The patch size ratio is $r=1/6$ and the mean slope of the bed is $\tan\theta=0.001$. 
}
\label{patchNp}
\end{figure}

For example, we present the case of \(2\pi\)-periodic dynamics  simulated by $m=10$ patches in a period, and $n=9$ interior microscale grid points on each patch.
The patch size ratio was \(r=1/6\) which means about a third of the spatial domain is covered by patches, and about two-thirds is unsimulated space: in practice we aim use smaller~\(r\) but this~\(r\) shows structures more clearly.
Setting the bed slope to \(\tan\theta=0.001\) an equilibrium flow  of equation~\eqref{patch:Nu} is that depth $h=1$ and depth-averaged lateral velocity $\uu\approx18.1\tan^{1/2}\theta=0.57$  (non-dimensional).
Figure~\ref{patchNp} plots a numerical gap-tooth simulation for the depth~$h$ and depth-averaged lateral velocity~$\uu$ at three times.
At the initial time $t=0$\,, we superimposed on the equilibrium flow a macroscale wave of~$0.2\sin x$ together with small random microscale noise. 
The $t=2$ graph shows that the microscale structures within a patch has smoothed quickly by the microscale dissipation~$h|\uu|\DD x\uu$.
%This corresponds to the large decay rates of the microscale modes in  Figure~\ref{patchNeig}.
In addition, the macroscale wave propagates downstream on the free surface, decaying slowly, as illustrated by the $t=4$ graph.
%This slow decay of macroscale waves corresponds to the eigenvalues of small real part but non-zero frequencies in Figure~\ref{patchNeig}.
%Section~\ref{sec:dam} presents and discusses a large length scale simulation, using the gap-tooth scheme with the nonlinear microscale turbulent model~\eqref{eqs:patch:Nji} and the coupling conditions~\eqref{patch:cph}.

\begin{figure}
\centering
\setlength{\unitlength}{0.7ex}
\begin{picture}(83,65)(10,-5)
%\put(10,-5){\framebox(83,65){}}
\put(15,0){\includegraphics[height=60\unitlength]{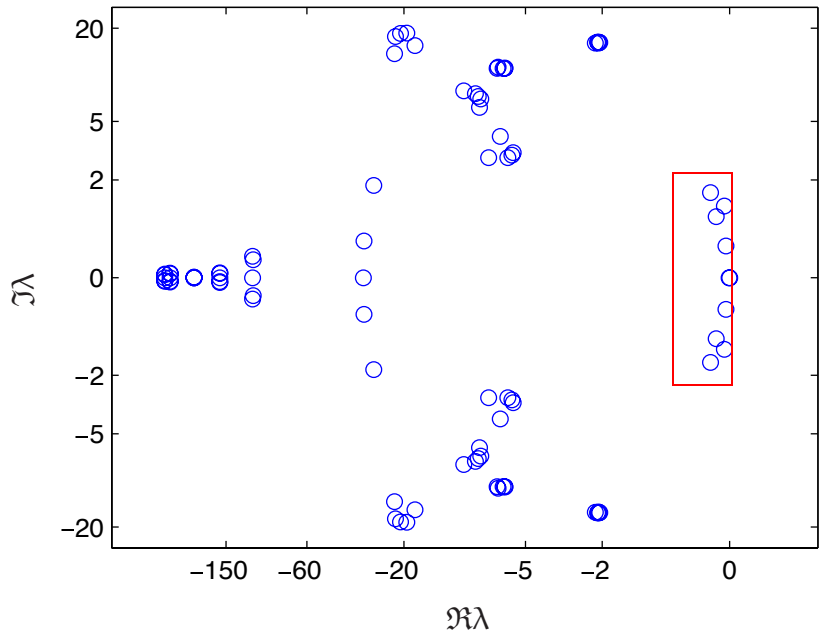}}
\put(10,0){\rotatebox{90}{\hspace{20ex}$\Im\lambda$}}
\put(55,-5){$\Re\lambda$}
\put(78.5,20.5){\color{red}\framebox(6,22)}
\end{picture}
\caption{Distribution of the real and imaginary parts of the numerical eigenvalues (non-uniform scaling of axes) for flood \pde{}s~\eqref{eqs:patch:Nji} with $m=10$ patches and $n=9$ microscale grid points on each patch. 
The domain~\(\Omega\) is assumed $2\pi$-periodic, the length scale ratio $r=1/6$\,, and the bed slope $\tan\theta=0.001$.
}
\label{patchNeig}
\end{figure}

\begin{figure}
\centering
\begin{tabular}{c@{\ }c}
\rotatebox{90}{\hspace{20ex}$\Im\lambda$} &
\includegraphics[height=45ex]{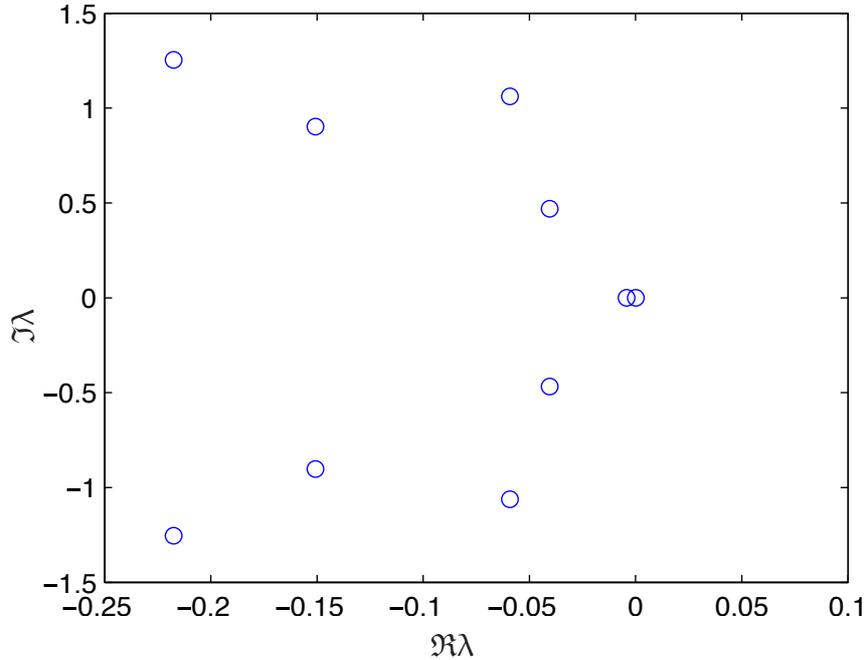}\\
& $\Re\lambda$
\end{tabular}
\caption{Zoom into the values in the red box in Figure~\ref{patchNeig}.
There are four complex conjugate pairs of values, and another two real values~$0.0002$ and $-0.0042$.
}
\label{patchNeigM}
\end{figure}

Eigenvalue analysis illuminates the gap-tooth scheme when applied to nonlinear water wave microscale discretisation~\eqref{eqs:patch:Nji} with the cubic coupling conditions~\eqref{eqs:cubic}.
%Consider the nonlinear microscale simulator~\eqref{eqs:patch:Nji} and the coupling conditions~\eqref{patch:cph} to the error~$\mathcal O(\gamma^3)$.  
%Omitting terms~$\mathcal O(\gamma^3)$, the boundary conditions~\eqref{patch:cph} become a cubic interpolation of the neighbouring macroscale values to approximate the values on a patch edges.
Let's consider further the waves upon the equilibria of water flowing down an inclined plane with bed drag balancing gravitational forcing.
By the non-dimensionalisation we just consider the equilibrium of water depth $h_{j,i}=1$ and depth-averaged lateral velocity $\uu_{j,i}\approx18.1\tan^{1/2}\theta$.
The spectrum of the Jacobian of the system characterises the dynamics in the neighbourhood of this equilibrium.
We estimated the Jacobian to about seven digits accuracy via centred numerical differentiation of the gap-tooth simulation procedure.  Then standard routines computed the complete spectrum of eigenvalues.

%\begin{table}
%\caption{Numerical eigenvalues for the macroscale wave modes from Figure~\ref{patchNeig} for \(m=10\) and \(m=14\) patches; compare with exact eigenvalues from equation~\eqref{patch:lambda} for the linear dynamics about the steady state.}
%\label{patchNeigTM}
%\begin{displaymath}
%\begin{array}{ccccc}
%\hline
%  k & \text{multiplicity} & 
%  \lambda, \text{ eqn~\eqref{patch:lambda}}& 
%  \lambda,\ m=10 &
%  \lambda,\ m=14\\
%\hline
%0 & \text{one}& 0,-0.003 &  0.0002 , -0.004& 0.0002, -0.005 \\
%\pm 1 & \text{one} &  -0.075\pm0.404i & -0.041 \pm 0.486i& -0.041\pm 0.502i \\
%\pm2 & \text{one} & -0.296\pm0.791i & -0.151 \pm 1.030i & -0.155\pm1.027i\\
%  \hline
%\end{array}
%\end{displaymath}
%\end{table}

For example, to match the simulation of Figure~\ref{patchNp} we analysed the case of \(2\pi\)-periodic dynamics  predicted by $m=10$ patches in a period, and $n=9$ interior microscale grid points on each patch to give a system in \({90}\)~variables.
The patch size ratio was \(r=1/6\)\,.
Figure~\ref{patchNeig} plots the growth rate~($\Re\lambda$) and frequency~($\Im\lambda$) obtained from the Jacobian.
There are $40$~pairs of complex conjugate eigenvalues with large negative real parts which represent the microscale modes within the patches.
Most of these negative real parts are between~$-2$ and~$-200$.
These represent microscale waves within the patches which decay rapidly through the microscale turbulent dissipation, dominantly~$h|\uu|\DD x\uu$.
The large imaginary parts of many of these eigenvalues reflect the fast oscillation of the microscale waves within the patches.

With \(m=10\) patches on a staggered macroscale grid the scheme resolves four macroscale waves: two propagating upstream and two downstream. 
The eigenvalues with small real parts in the red box of Figure~\ref{patchNeig} represent these macroscale waves.
Figure~\ref{patchNeigM} zooms in these small eigenvalues and shows that there are four pairs of small decay, the different imaginary-part (frequencies) distinguish upstream and downstream propagation, together with two real values near zero ($0.0002$ and~$-0.0042$).
This slow decay of the four macroscale waves is due to a combination of the turbulent dissipation and the small nonlinear bed drag~$0.003\uu|\uu|/h$.
The small negative eigenvalue of~\(-0.0042\) represents the decay of uniform flow to the equilibrium speed down the sloping bed.
The smallest eigenvalue is zero to numerical error in the Jacobian and represents conservation of water.
This spectrum of eigenvalues confirms the gap-tooth scheme models the macroscale wave propagation without the scheme explicitly knowing any macroscale closure.

The gap between the growth rate~$\Re\lambda\approx0$ and~$\Re\lambda\approx-2$, as shown in Figure~\ref{patchNeig}, characterises the separation between slow macroscale dynamics and the fast microscale dynamics.
As discussed in section~\ref{sec:nsm}, this gap indicates that there is a nonlinear slow manifold of the macroscale modes in the gap-tooth dynamics \cite[e.g.]{Roberts1988, Chicone:2006fk, Potzsche:2006uq}.
Sections~\ref{sec:consistency}--\ref{sec:consistencyC} first establish the high order consistency between the macroscale gap-tooth scheme and the underlying microscale system.

\section{Coupling conditions connect patches across space}
\label{sec:couple}

We propose and analyse a macroscale staggered grid of patches as shown in Figure~\ref{patchMacro}.
This section establishes that classic macroscale interpolation provides a coupling between patches that achieves arbitrarily high order consistency between the patch scheme and the underlying microscale dynamics for general microscale systems.

\subsection{Prove consistency for general wave systems}
\label{sec:consistency}

%??
%However, the coupling between patches is best defined in terms of the microscale fields to form a closed system across the whole domain.
%??

For definiteness in theoretical support, let there be \(m\)~patches in a spatial domain~\(\Omega\) with the fields required to be periodic in~\(x\), and the fields to be in the Sobolev space~\(\HH^2(\Omega)\) of square-integrable functions on~\(\Omega\) with square-integrable first two derivatives.
We model the dynamics on patches~\(E_j\) of an equi-spaced macroscale grid $X_j=jD$\,: define the patches \(E_j:=\{x\in\Omega \mid |x-X_j|<rD \}\), for \(j=1,\ldots,m\), centred on each macroscale grid point, and the collection of patches \(E:=\{E_j \mid j=1,\ldots,m\}\).
The parameter~$r$ is the ratio between a patch half-width and the macroscale step~$D$.
The parameter~\(r\) characterises the size of each patch relative to the distance between neighbouring patches: 
when $r=1/2$ the neighbouring patches meet; 
and when $r=1$ the patches overlap as was found so useful in holistic discretisation \cite{Roberts1999, Roberts00a}.
When the ratio~$r$ is small, the patches form a relatively small part of the physical domain to engender a computationally efficient scheme for multiscale simulation.

To couple patches of microscale simuations, define a macroscale shift operator \(\shift h(x):=h(x+2D)\)
\footnote{The full shift~\shift\ skips a patch because we employ a staggered macroscale grid as shown in Figure~\ref{patchMacro}.  Expressions for mean and difference operators have to be correspondingly adjusted for the staggered grid.} and corresponding centred difference and mean operators here defined as \(\delta:=\shift^{1/2}-\shift^{-1/2}\) and \(\mu:=(\shift^{1/2}+\shift^{-1/2})/2\).
Consider a field~\(h(x)\in\HH^2(\Omega)\).  
Evaluating the field~\(h\) at a shift corresponding to the width of a patch  gives
\begin{align}
h(x\pm rD)&=\shift^{\pm r/2}h
\nonumber\\&
=\left(1\pm\mu\delta+\tfrac12\delta^2\right)^{r/2}h
\quad\text{\cite[p.65, e.g.]{npl61}}
\nonumber\\&
=\frac{\mu}{\sqrt{1+\delta^2/4}}(1\pm\mu\delta+\tfrac12\delta^2)^{r/2}h
\quad(\text{as }\mu^2=1+\delta^2/4)
\nonumber\\&
=\left[\mu\pm\tfrac12r\delta+\tfrac18(-1+r^2)\mu\delta^2\pm\tfrac{1}{48}(-r+r^3)\delta^3+\cdots\right]h\,,\label{eq:shif}
\end{align}
This expansion motivates defining the corresponding ameliorated coupling operators
\begin{align}
\cC_{\pm}:={}&
\gamma\left[\mu\pm\tfrac12 r\delta\right]+\gamma^3\left[\tfrac18(-1+r^2)\mu\delta^2\pm\tfrac{1}{48}(-r+r^3)\delta^3
\right]
\nonumber\\&{}
+\gamma^5\left[\tfrac{1}{384}(9-10r^2+r^4)\mu\delta^4\pm\tfrac{1}{3840}(9r-10r^3+r^5)\delta^5\right]
\nonumber\\&{}
+\gamma^7\left[\tfrac{1}{46080}(-225+259r^2-35r^4+r^6)\mu\delta^6
\right.\nonumber\\&\quad\left.{}
\pm\tfrac{1}{645120}(-225r+259r^3-35r^5+r^7)\delta^7\right]
+\cdots\,,\label{patch:cphco}
\end{align}
that operates on the patch index~\(j\), in terms of mean and difference operators that hereafter operate on index~\(j\): for example, \(\delta H_j=H_{j+1}-H_{j-1}\) and \(\mu H_j=(H_{j+1}+H_{j-1})/2\).
The parameter~$\gamma$ conveniently labels the spatial extent of the various terms appearing in coupling operators~\eqref{patch:cphco} (discussed in the next paragraph).
Then, as suggested by Figure~\ref{patchMacro}, we require boundary conditions for each patch of
\begin{align}
\uu_j(X_j\pm rD,t)=\cC_\pm \uu_j(X_j,t)
\quad\text{and}\quad
h_j(X_j\pm rD,t)=\cC_\pm h_j(X_j,t),
\label{patch:cph}
\end{align}
for odd and even~\(j\) respectively.
By the definition of macroscale values~\eqref{patch:macrogrid}, the right-hand sides of the boundary conditions~\eqref{patch:cph} couple the microscale patches together via interpolation of the macroscale grid values.

We obtain various accuracies for the macroscale simulation by truncating the coupling~\eqref{patch:cph} to various orders in the label~$\gamma$.
This follows as \(\gamma\)~parametrises the stencil width of the interpolation.
For example, truncating to errors~\Ord{\gamma^3} gives linear interpolation from the nearest neighbour patches, \(h_j(X_j\pm rD,t)=(\mu\pm\tfrac12r\delta)H_j\), as Figure~\ref{patchMacro} illustrates.
Whereas truncating to errors~\Ord{\gamma^5} gives cubic interpolation from nearest and next nearest patches.
Truncating to errors~\Ord{\gamma^7} gives a quintic interpolation, and so on.
One key property of the coupling is that, when truncated to errors~\Ord{\gamma^p}, and upon setting label \(\gamma=1\)\,, the difference \(\cC_+-\cC_-=\shift^{r/2}-\shift^{-r/2}+\Ord{\delta^{p-1}}\) in the limit of large length scale macroscale variations.
Next we prove that the various order interpolations of the gap-tooth scheme achieves corresponding orders of consistency with the microscale system over~\(\Omega\)---a consistency analogous to that for dissipative systems \cite[e.g.]{Roberts:2005uq, Roberts:2011uq}.

\begin{theorem} \label{thm:consistency}
Consider the general coupled system of two equations
\begin{equation}
\partial_th=\cL_1\uu+c_1h%\label{eq:Leqhu}
\quad\text{and}\quad
\partial_t\uu=\cL_2h+c_2\uu\,,\label{eq:Leqhu}
\end{equation}
for fields \(h(x,t),\uu(x,t)\in\HH^2(\Omega)\),
for generic odd, homogeneous linear operators~$\cL_1$ and~$\cL_2$, and constants~$c_1$ and~$c_2$.
Let $h_j(x,t),\uu_j(x,t)\in\HH^2(E_j)$ denote the subgrid fields on the $j$th~patch satisfying~\eqref{eq:Leqhu} in patches~$E_j$ with the coupling conditions~\eqref{patch:cph}.
When inter-patch coupling~\eqref{patch:cph} is truncated to residuals~\Ord{\gamma^p}, then the macroscale grid values~\eqref{patch:macrogrid} evolve consistently with the equations~\eqref{eq:Leqhu} to errors~\Ord{\delta^{p-1}} (upon setting \(\gamma=1\)).
\end{theorem}

In this theorem the linear operators~\(\cL_\ell\) could represent partial derivatives, microscale discretisations, lattice Boltzmann interactions, and so on.
The following proof is based upon the ideas used to prove analogous consistency for dissipative systems \cite[Theorem~7]{Roberts:2011uq}.

\begin{proof} 
The proofs for each of the equations~\eqref{eq:Leqhu} are the same with appropriate interchange of symbols.
Consider \(\partial_t\uu=\cL_2h+c_2\uu\)\,.
Because the linear operator~$\cL_2$ is odd and homogeneous, we formally expand the operator
\begin{equation}
\cL_2=\sum_{k=0}^{\infty}{\ell_{2k+1}\delta_r^{2k+1}}=\ell(\delta_r)\,,
\label{eq:linearL}
\end{equation}
in terms of the patch sized, microscale, centred difference \(\delta_r:=\shift^{r/2}-\shift^{-r/2}\),
for some coefficients~$\ell_{2k+1}$ and corresponding function~$\ell$.
By the term `generic odd' in the theorem, we mean the coefficient \(\ell_1\neq 0\)\,.
The second of~\eqref{eq:Leqhu} on the patch~\(E_j\) determines \(\partial_t\uu_j=\cL_2h_j+c_2\uu_j\) and now becomes $(\partial_t-c_2)\uu_j=\ell(\delta_r)h_j$\,.
Because \(\ell_1\neq 0\)\,, function~\(\ell\) has a smooth inverse function~\(\ell^{-1}\), at least near zero, and so we rearrange this microscale equation to
\begin{align}
\ell^{-1}(\partial_t-c_2)\uu_j=\delta_rh_j\,.
\label{eq:Linverse}
\end{align}
Now evaluate~\eqref{eq:Linverse} at the patch centre $x=X_j$\,: on the left-hand side the time derivatives commute with the evaluation at \({x=X_j}\) so equation~\eqref{eq:Linverse} becomes, by definition~\eqref{patch:macrogrid},
\begin{align}
\ell^{-1}(\partial_t-c_2)U_j
&=h_j(X_j+rD,t)-h_j(X_j-rD,t)
\nonumber\\&
=\cC^+h_j(X_j,t)-\cC^-h_j(X_j,t)
\quad(\text{by coupling~\eqref{patch:cph}})
\nonumber\\&
=(\cC^+-\cC^-)H_j
\quad(\text{by definition~\eqref{patch:macrogrid}})
\nonumber\\&
=(\shift^{+r/2}-\shift^{-r/2})H_j+\Ord{\delta^{p-1}H_j}
\quad(\text{by truncating~\eqref{patch:cphco} with }\gamma=1)
\nonumber\\&
=\delta_rH_j+\Ord{\delta^{p-1}H_j},
\label{eq:LinearC}
\end{align}
where~$p$ is the order of error in~\(\gamma\) of the truncated coupling operators~\eqref{patch:cphco} (if order~\(p\) is even then the error is~\Ord{\delta^{p}H_j}).
Equation~\eqref{eq:LinearC} is a closed relation among the macroscale quantities.
Reverting the inverse function~$\ell^{-1}$, and equation~\eqref{eq:LinearC} implies
\begin{align*}
(\partial_t-c_2)U_j&=\ell(\delta_r)H_j+\Ord{\delta^{p-1}H_j}=\cL_2H_j+\Ord{\delta^{p-1}H_j},
\end{align*}
which then becomes
\begin{align}
\partial_tU_j=\cL_2H_j+c_2U_j+\Ord{\delta^{p-1}H_j}.
\label{eq:consiscp}
\end{align}
Similarly for the companion equation of~\eqref{eq:Leqhu}.
That is,  in the patch scheme with coupling conditions~\eqref{patch:cph}, the macroscale grid values~\eqref{patch:macrogrid} evolve consistently to any specified order with the microscale system~\eqref{eq:Leqhu} solved on the whole domain~\(\Omega\).
\end{proof}

\subsection{Computer algebra establishes further consistency}
\label{sec:consistencyC}

The previous subsection established consistency for general linear wave systems with simple drag, whereas we generally want to apply the patch scheme to wave systems with other dissipative mechanisms, and to nonlinear systems.
This section uses computer algebra to show that consistency is also obtained for a variety of such interesting systems.

\subsubsection{Algebraically confirm Theorem~\ref{thm:consistency}}
\label{sec:act1}

Consider the following dispersive system in the wave-like form~\eqref{eq:Leqhu}:
\begin{align}&
\D th=c_1h-\D x\uu-c_{11}\DDD{x}{\uu}
%\label{eq:case1h}\\&
\quad\text{and}\quad
\D t\uu=c_2\uu-\D xh-c_{21}\DDD xh\,,
\label{eq:case1}
\end{align}
with  constant coefficients $c_1,c_{11},c_2,c_{21}$.
We confirm the consistency, established by Theorem~\ref{thm:consistency}, between the gap-tooth scheme and this underlying microscale system.
%Centre manifold theorems~\cite[e.g.]{Roberts1988,Potzsche:2006uq} demonstrate that there is a slow manifold in this multiscale system.
Computer algebra constructs solutions to the system~\eqref{eq:case1} on the patches~\(E_j\) when coupled by~\eqref{patch:cph} as a regular power series in the coupling parameter~\(\gamma\).
In the solution, terms of up to~\(\gamma^{p-1}\) then encode all the effects of truncating the coupling condition~\eqref{patch:cph} to errors~\Ord{\gamma^p}; that is, to account for interactions between a patch and its \(p\)~neighbours on either side.

For example, the computer algebra of Appendix~\ref{sec:appendix} (with \verb|choice:=1|) derives the microscale field in each patch as \((h_j,\uu_j)=\cM(H_j,U_j)\) where, in terms of \(\xi=(x-X_j)/(rD)\)\,, operator
\begin{equation*}
\cM=\begin{cases}
1
+\gamma^2\left[\tfrac12r\xi\delta\mu+\tfrac18r^2\xi^2\delta^2\right]
+\Ord{\gamma^4},
\\
\gamma\left[\mu+\tfrac12r\xi\delta\right]
+\gamma^3\left[\tfrac18(-1+r^2\xi^2)\delta^2\mu+\tfrac1{48}(-r\xi+r^3\xi^3)\delta^3\right]
+\Ord{\gamma^4},
\end{cases}
\end{equation*}
alternating upon whether \(j\)~is even or odd and whether applied to \(h_j\) or~\(\uu_j\).
Remarkably, in the class of \pde{}s~\eqref{eq:case1} the microscale field is independent of the coefficients~\(c_k\): this independence does not generally occur in the other classes of \pde{}s.
For these microscale fields, the computer algebra, Appendix~\ref{sec:appendix}, derives the corresponding evolution of the macroscale values to be, for the appropriate~\(k\),
\begin{equation}
(\dot H_j,\dot U_j)=\left[
-\frac1D(\tfrac12\gamma\delta-\tfrac1{48}\gamma^3\delta^3)
-\frac1{8D^3}c_{k1}\gamma^3\delta^3\right](U_j,H_j)
+c_k(H_j,U_j)+\Ord{\gamma^4},
\label{eq:macroeasy}
\end{equation}
These macroscale evolution equations correspond to a conventional macroscale discretisation of the microscale \pde{}s~\eqref{eq:case1}.
But remember that the gap-tooth scheme would generate a macroscale simulation obeying~\eqref{eq:macroeasy} without knowing explicitly such a closure.

The required high order consistency to confirm Theorem~\ref{thm:consistency} is explored by transforming such macroscale discrete models~\eqref{eq:macroeasy} to its equivalent \pde\ and comparing to the microscale \pde~\eqref{eq:case1}.
Truncating the coupling conditions~\eqref{patch:cph} to errors~\Ord{\gamma^9}, in a couple of \textsc{cpu} seconds the computer algebra program of Appendix~\ref{sec:appendix} derives the higher order version of the macroscale model~\eqref{eq:macroeasy}.
Post-processing then uses Taylor series, \(H_{j+p}=\sum_{n=0}^\infty (pD)^n/n!\,\Dn XnH\)\,, to transform the higher order version of~\eqref{eq:macroeasy} to the equivalent \pde{}s for the macroscale variables as
\begin{eqnarray}
\D tH&=&-\left[\gamma\D XU
+\tfrac16(\gamma-\gamma^3)D^2\Dn X3U
+\tfrac1{120}(\gamma-10\gamma^3+9\gamma^5)D^4\Dn X5U
\right.\nonumber\\&&\quad\left.{}
+\tfrac1{5040}(\gamma-91\gamma^3+315\gamma^5-225\gamma^7)D^6\Dn X7U
\right]
\nonumber\\&&{}
-c_{11}\left[\gamma^3\Dn X3U
+\tfrac12(\gamma^3-\gamma^5)D^2\Dn X5U
+\tfrac1{120}(13\gamma^3-50\gamma^5+37\gamma^7)D^4\Dn X7U
\right]
\nonumber\\&&{}
+c_1H+\Ord{(D^8+c_{11}D^6)\Dn X9U},
\label{eq:conthmh}
\end{eqnarray}
and similarly for \(\D tU\).
Observe in these expressions how beautifully various contributions cancel when artificial parameter~\(\gamma\) is set to one: the result is that the macroscale variables~\(H\) and~\(U\) evolve consistently with the microscale \pde{}s~\eqref{eq:case1}. 
However, when the coupling between patches is limited to \(p-1\)~nearest neighbouring patches on either side, equivalent to truncating the coupling~\eqref{patch:cphco} to errors~\Ord{\gamma^p}, then the \pde~\eqref{eq:conthmh} confirms the consistency holds to errors~\Ord{\Dn XpU} for \(p\in\{3,5,7,9\}\); that is, to errors~\Ord{\delta^pU_j} in accord with Theorem~\ref{thm:consistency}.

\subsubsection{Linear waves with dissipation}

Consider the following coupled system supporting dispersive waves that are damped by a diffusion of strength \(c_3\) and~\(c_4\):
\begin{equation}
\D th=-\D{x}{\uu}-c_{11}\DDD{x}{\uu}+c_3\DD xh
%\label{eq:case2h}
\quad\text{and}\quad
\D{t}{\uu}=-\D{x}{h}-c_{21}\DDD{x}{h}+c_4\DD{x}{\uu}\,.
\label{eq:case2}
\end{equation}
Because of the diffusive dissipation, this system is not in the form~\eqref{eq:Leqhu} addressed by Theorem~\ref{thm:consistency}, yet we here demonstrate similar high order consistency between a gap-tooth scheme for this system and the \pde{}s.

The computer algebra program of Appendix~\ref{sec:appendix} (with \verb|choice:=2|) solves the \pde{}s~\eqref{eq:case2} on patches coupled by~\eqref{patch:cphco}.
The resultant microscale fields in each patch and the discrete evolution are analogous to that obtained in section~\ref{sec:act1}.
Thus the equivalent \pde\ of the macroscale evolution is also the same as~\eqref{eq:conthmh} except for additional terms introduced by the diffusive dissipation:
\begin{eqnarray}
\D tH&=&(\text{first three lines of~\eqref{eq:conthmh}})
%-\left[\gamma\D XU
%+\tfrac16(\gamma-\gamma^3)D^2\Dn X3U
%+\tfrac1{120}(\gamma-10\gamma^3+9\gamma^5)D^4\Dn X5U
%\right.\nonumber\\&&\quad\left.{}
%+\tfrac1{5040}(\gamma-91\gamma^3+315\gamma^5-225\gamma^7)D^6\Dn X7U
%\right]
%\nonumber\\&&{}
%-c_{11}\left[\gamma^3\Dn X3U
%+\tfrac12(\gamma^3-\gamma^5)D^2\Dn X5U
%+\tfrac1{120}(13\gamma^3-50\gamma^5+37\gamma^7)D^4\Dn X7U
%\right]
\nonumber\\&&{}
+c_3\left[\gamma^2\DD XH
+\tfrac13(\gamma^2-\gamma^4)D^2\Dn X4H
+\tfrac2{45}(\gamma^2-5\gamma^4+4\gamma^4)D^4\Dn X6H
\right.\nonumber\\&&\left.\quad{}
+\tfrac1{315}(\gamma^2-21\gamma^4+56\gamma^6-36\gamma^8)D^6\Dn X8H
\right]
\nonumber\\&&{}
+\Ord{(D^8+c_{11}D^6)\Dn X9U+c_3D^8\Dn X{10}H},
\label{eq:case2hf}
\end{eqnarray}
and similarly for \(\D tU\).
As before, when the coupling between patches is limited to \(p-1\)~nearest neighbouring patches on either side, equivalent to truncating the coupling~\eqref{patch:cphco} to errors~\Ord{\gamma^p}, then the \pde~\eqref{eq:case2hf} confirms that consistency holds to errors~\Ord{\Dn Xp{}} for \(p\in\{3,5,7,9\}\); that is, to errors~\Ord{\delta^p}.
This establishes the consistency of the gap-tooth scheme for the class of dispersive wave systems~\eqref{eq:case2} with diffusive dissipation.

\subsubsection{Nonlinear wave systems}

We generally want to use the gap-tooth scheme for macroscale simulation of nonlinear microscale dynamics.
This section explores the basic example nonlinear wave system
\begin{equation}
\D th=-\D{x}{\uu}
\quad\text{and}\quad
\D{t}{\uu}=-\D xh-c_5\uu\D{x}{\uu}\,,
\label{eq:case3hu}
\end{equation}
where \(\uu\D x\uu\) in the second \pde\ is typical of the self-advection of momentum.
Because of the nonlinearity, this system is not in the form~\eqref{eq:Leqhu} addressed by Theorem~\ref{thm:consistency}, yet again we here demonstrate consistency between a gap-tooth scheme for this system and the \pde{}s.

The computer algebra program of Appendix~\ref{sec:appendix} (with \verb|choice:=3|) solves the nonlinear \pde{}s~\eqref{eq:case3hu} on patches coupled by~\eqref{patch:cphco}.
But the nonlinearity needs to be small, and the algebraic complexity is great, so the computer algebra here solves the \pde{}s~\eqref{eq:case3hu} to errors~\Ord{\gamma^7,c_5^3}.
Expressions for the microscale field and evolution of the macroscale quantities are rather complicated, and so omitted.
Instead, we just record part of the equivalent \pde{}s for the macroscale variables in order to indicate how consistency with the microscale develops.
Limiting the expressions to errors~\Ord{D^4},
\begin{eqnarray}
H_t&=&-\gamma U_X
+\tfrac16D^2\gamma(1-\gamma)\left[\vphantom{H_X^1}
-(1+\gamma)U_{XXX}
\right.\nonumber\\&&\left.\quad{}
+c_5\gamma^2r^2\big(-2H_{XX}U_X+\gamma \{H_XU_{XX}+ H_{XXX}U\}\big)
+c_5^2\gamma^3r^2\big(-2U_X^3
\right.\nonumber\\&&\left.\qquad{}
-2U_XU_{XX}
+\gamma\{5U_{XX}U_X+2U^2U_{XXX}+H_X^2U_X+3H_{XX}H_XU\}\big)
\right]
\nonumber\\&&{}
+\Ord{\gamma^7,c_5^3,D^4},
\label{eq:hhtn}
\\U_t&=& -\gamma H_X -\gamma^2 c_5UU_X
+\tfrac16D^2\gamma(1-\gamma)\left[\vphantom{H_X^1}
(1+\gamma)H_{XXX}
\right.\nonumber\\&&\left.\qquad{}
+c_5^2\gamma^2r^2\big(
+2\gamma H_{XX}UU_X
+\{2+2\gamma-\gamma^2\}\{H_XUU_{XX}+H_{XXX}U^2\}
\big)
\right]
\nonumber\\&&{}
-c_5D^2\gamma^2\tfrac1{24}\big(16-16\gamma^2+\gamma^4\{9-10r^2+r^4\}\big)UU_{XXX}
\nonumber\\&&{}
+\Ord{\gamma^7,c_5^3,D^4}.
\label{eq:uutn}
\end{eqnarray}
These equivalent \pde{}s for the macroscale dynamics of the gap-tooth scheme have several interesting aspects.
The components \(H_t=-\gamma U_X\) and \(U_t=-\gamma H_X -\gamma^2 c_5UU_X\) show that provided the inter-patch coupling involves at least the next nearest neighbours (flagged by the \(\gamma^2\)-factor), then the macroscale dynamics of the gap-tooth scheme is consistent with the microscale~\eqref{eq:case3hu} to errors indicated by the other terms.
The error terms of the \(H_t\)-\pde~\eqref{eq:hhtn} all have a factor~\((1-\gamma)\) so that these errors vanish when the inter-patch coupling is carried out to high enough order:  consequently consistency follows to~\Ord{D^4}.
Similarly for the first group of error terms of the \(U_t\)-\pde~\eqref{eq:uutn}.
However, the second group of error terms of~\eqref{eq:uutn}, that in~\(c_5D^2\), do not appear to have a factor~\((1-\gamma)\), and so need not vanish: potentially, higher order analysis could find the requisite factor and remove the error; alternatively, the error vanishes for the case of cubic interpolation when truncating inter-patch coupling~\eqref{patch:cphco} to~\Ord{\gamma^5}.
Interestingly, this consistency error also vanishes for the overlapping patch case \(r=1\) that is so attractive for `holistic discretisation' \cite[e.g]{Roberts2011a}, but overlapping patches are not relevant for efficient numerical simulation using the gap-tooth scheme.
Summarising simply, the macroscale dynamics of the gap-tooth scheme is consistent with the microscale nonlinear dynamics of~\eqref{eq:case3hu} to errors~\Ord{D^2}, at most.

\subsection{Nonlinear slow manifolds exist}
\label{sec:nsm}

The computer algebra of section~\ref{sec:consistencyC} constructs the macroscale dynamics as a slow manifold of the gap-tooth scheme \cite[e.g.]{Boyd95, MacKay03}.
In the case of linear systems the slow manifolds are more specifically slow subspaces.
This section establishes that such slow manifolds exist for some system close to that specified \cite[Chapt.~13]{Roberts2014a}, and identifies that the fast microscale waves, if undamped, may nonetheless affect the macroscale dynamics.

We establish a gap-tooth slow manifold for microscale systems in the general form
\begin{equation}
\partial_th=\cL_1\uu+\epsilon f_1(h,\uu)
\quad\text{and}\quad
\partial_t\uu=\cL_2h+\epsilon f_2(h,\uu)\,,\label{eq:hunls}
\end{equation}
for fields \(h(x,t),\uu(x,t)\in\HH^2(\Omega)\), some sufficiently smooth functions~\(f_\ell\), and homogeneous operators~\(\cL_\ell\) satisfying the following fast-slow dichotomy:
forming the operators~\(\cL_\ell\), whether differential as in~\eqref{eqs:patch:N} or discrete as in~\eqref{eqs:patch:Nji}, 
into the two combined operators
\begin{equation*}
\cK_\ell=\begin{bmatrix} 0&\cL_\ell\\\cL_{3-\ell}&0 \end{bmatrix}
\end{equation*}
such that the second component satisfies Dirichlet boundary conditions at edges of some interval, then both \(\cK_1,\cK_2\) must have a zero eigenvalue with corresponding eigenvector~\((1,0)\), and all other eigenvalues~\(\lambda\) must be bounded away from zero, in modulus, \(|\lambda|>\beta>0\)\,.
Form a gap-tooth system by letting $h_j(x,t),\uu_j(x,t)\in\HH^2(E_j)$ denote the subgrid fields on the $j$th~patch satisfying~\eqref{eq:hunls} in patches~$E_j$ with the coupling conditions~\eqref{patch:cph}.
Because of the zero eigenvalue of~\(\cK_\ell\), for no nonlinearity, \(\epsilon=0\)\,, and no coupling, \(\gamma=0\)\,, this gap-tooth system has a subspace of equilibria~\(\MM_0\) of piecewise constant fields in each patch:
\begin{equation*}
(h_j,\uu_j)=(H_j,0)\text{ for odd }j,\quad
(h_j,\uu_j)=(0,U_j)\text{ for even }j,
\end{equation*}
for independent `amplitudes'~\(H_j\) and~\(U_j\).
Set the eigenvalue bound~\(\beta\) to be the smallest necessary for the range of~\(\MM_0\) of interest; typically \(\beta\propto1/r\)\,.
Then a power series construction finds the following \cite[\S3]{Cox93b}, \cite[Chapt.~13]{Roberts2014a}: based at each of these equilibria there exists a smooth system and a smooth coordinate transformation which, firstly, together are \Ord{\gamma^p,\epsilon^q}-close to the gap-tooth system~\eqref{eq:hunls} with coupling~\eqref{patch:cph}, and secondly, possesses a slow manifold~\cM\ global in~\((\vec H,\vec U)\).
\footnote{Notice that this statement is a `backwards theory' \cite[e.g.]{Grcar2011} that neatly sidesteps the controversy over the existence or otherwise of slow manifolds \cite[e.g.]{Lorenz86, Lorenz87, Lorenz92, Jacobs91}.}
Section~\ref{sec:consistencyC} constructed and discussed that part of the coordinate transform that was on the slow manifold itself.
This theory asserts that the macroscale grid values~\(H_j\) and~\(U_j\) discussed in section~\ref{sec:consistencyC} can form a sound closure for a finite range of nonlinearity~\(\epsilon\) and coupling~\(\gamma\).

However, for wave systems there is no assurance that the slow manifold is emergent. 
For nonlinear wave-like systems, the long term evolution on and off the slow manifold~\cM\ may be different---generally different by an amount quadratic in the fast waves \cite[]{Cox93b}.
A user needs to be wary if the fast microscale waves in the patches~\(E\) persist as a significant feature of the dynamics.
The possibility is that then such fast waves `trapped' in~\(E\), through resonance, may affect the macroscale evolution in a way significantly different from the way such fast waves affect the evolution if distributed over all space~\(\Omega\). 
In our simulations of turbulent floods this is not a problem as the turbulent dissipation, \(h|\uu|\uu_{xx}\), that is so weak as to be typically negligible over interesting macroscales, is reasonably strong inside the microscale patches and damps the microscale fast waves.
Similarly, in many applications some dissipation that is negligible on the macroscale will be a significant dissipation on the microscale and so damps the fast waves to leave the system on the slow manifold~\cM.

\section{Gap-tooth simulation of dam breaking }
\label{sec:dam}

This section applies the gap-tooth scheme to simulate dam-breaking.
The aim is to show how the scheme caters for the extreme nonlinearity of the dam and turbulent bore discontinuities (section~\ref{patch:damN}) and to discuss ways to implement practical domain boundary conditions (section~\ref{patch:damB}) rather than the periodic domain used in previous Sections~\ref{sec:Eig} and~\ref{sec:couple}.

\begin{figure}
\begin{center}
\setlength{\unitlength}{0.7ex}
\begin{picture}(105,60)(-5,0)
%\put(-5,0){\framebox(105,60){}}
\put(0,0){\includegraphics[height=40ex]{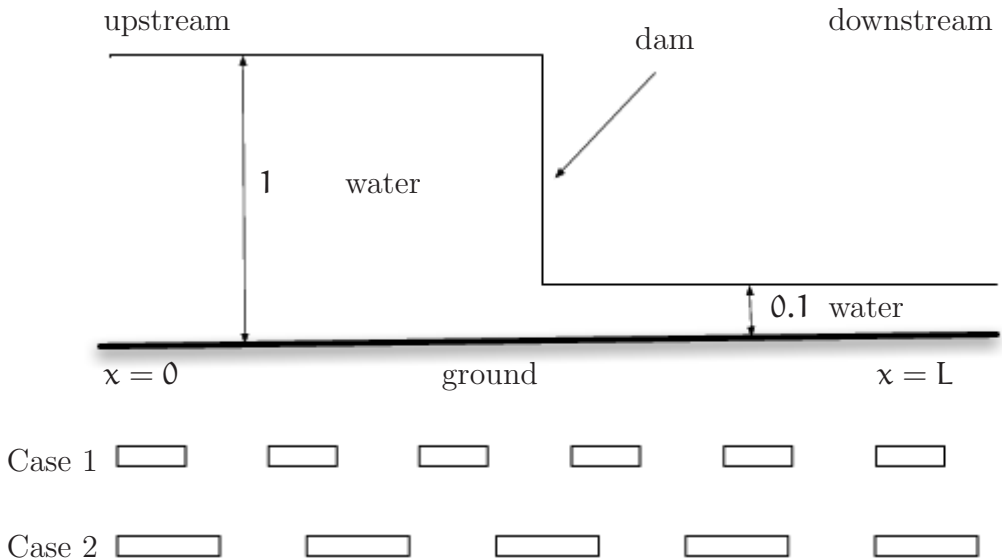}}
\put(21,40){$1$}
\put(74,27){$0.1$}
\put(5,20){$x=0$}
\put(85,20){$x=L$}
\put(30,40){{water}}
\put(80,27){{water}}
\put(5,57){{upstream}}
\put(80,57){{downstream}}
\put(60,55){{dam}}
\put(40,20){{ground}}
\put(-5,11){{Case~1}}
\put(-5,2){{Case~2}}
\end{picture}
\end{center}
\caption{The initial conditions of the dam breaking on a domain of nondimensional length~$L$.
The dam located at $x=L/2$ holds water of nondimensional depth $h=1$ in the upstream and a nondimensional shallow depth $h=0.1$ in the downstream.
Case~1 places the dam lies in between patches, and Case~2 places the dam in the middle of a patch.}
\label{patchDam}
\end{figure}

Figure~\ref{patchDam} shows a dam that holds back water standing in the middle of the domain of length~$L$. 
The ground is horizontal and let~$x$ denote the horizontal position.
Initially the dam holds water upstream of nondimensional depth $h=1$.
To avoid poor conditioning in the numerical calculation, downstream of the dam let the water have a shallow depth, for example~$0.1$.
At time $t = 0$ the dam breaks and the upstream water rushes downstream.

We simulate the dam-breaking waves by the gap-tooth scheme with the microscale turbulent model~\eqref{eqs:patch:N} and the cubic coupling conditions~\eqref{eqs:cubic}.
For comparison with the gap-tooth simulations, we also compute the microscale simulation over the whole domain, and  report experimental data from \cite{Stansby1998}
This section shows that the gap-tooth scheme reasonably simulates the dam-breaking waves.

%\subsection{Distribute patches of the dam breaking}
%\label{patch:damP}

There are typically two ways to distribute patches in the macroscale domain: either a patch includes the dam and the microscale resolves the sharp change in water depth, or the dam lies between two patches and the depth change is resolved only in the macroscale interpolation.
Case~1 of Figure~\ref{patchDam} distributes six patches in the macroscale domain.
The dam stands in the middle of the gap between the third and fourth patches.
An advantage of such a choice is to avoid the sharp discontinuity at the dam, because it arises in the gap which is not represented in the gap-tooth simulation.
Case~2 in Figure~\ref{patchDam} distributes five patches in the whole domain.
The dam is included on the centre of the third patch.
This choice would resolve the microscale details of the dynamics at the dam when the dam breaks.

\subsection{Numerical gap-tooth simulation of dam breaking}
\label{patch:damN}

This section explores numerical gap-tooth simulations of the dam-breaking waves.
\cite{Georgiev2008} used a previous version of the turbulent water wave \pde{}s~\eqref{eqs:patch:N} in simulating dam breaking; but they simulated the system over all space in the domain, not by the gap-tooth scheme.
The comparison of the calculations and experimental data by~\cite{Georgiev2008} shows that the \pde{}s~\eqref{eqs:patch:N} are a reasonable  model of dam breaking waves.
This section compares the gap-tooth simulation with the microscale simulation over the whole space domain, and with some experimental data of \cite{Stansby1998}.
The simulations indicate that putting the dam within a patch appears better.

We compute both the gap-tooth simulation and the microscale simulation over the whole domain for the dam-breaking waves.
In the experiments by~\cite{Stansby1998}, the initial water depth behind the dam is~$10\cm$ and lies in a horizontal domain of length of~$200\cm$. 
The dam stands at the centre of the domain.
For comparison, we nondimensionalise the depth~$10\cm$ to one; then the nondimensional length is $L=200\cm/10\cm=20$.
In this pilot study of gap-tooth simulation, we distributed both $m=10$ and $m=22$ patches on the whole macroscale domain and use $n=9$ microscale grid points on a patch: \(m=10\)~patches is low resolution of the macroscale, and \(m=22\)~patches is only a medium resolustion.
Then the distance between neighbouring patches is $D=L/m=2$.
For the low resolution simulation, using the scale ratio $r=1/6$,  the width of each patch is $\ell=2rD=0.67$\,, and the microscale grid step in each patch is $d=\ell/(n+1)\approx0.07$.
For consistent comparison, let the spatial step in the microscale simulation over the whole macroscale space domain (not by the gap-tooth scheme) have the same microscale spatial step~$d$.

\begin{figure}
\centering
\begin{tabular}{r@{\ }c}
$t=0.0$\\
\rotatebox{90}{\hspace{10ex}$h$} &
\includegraphics{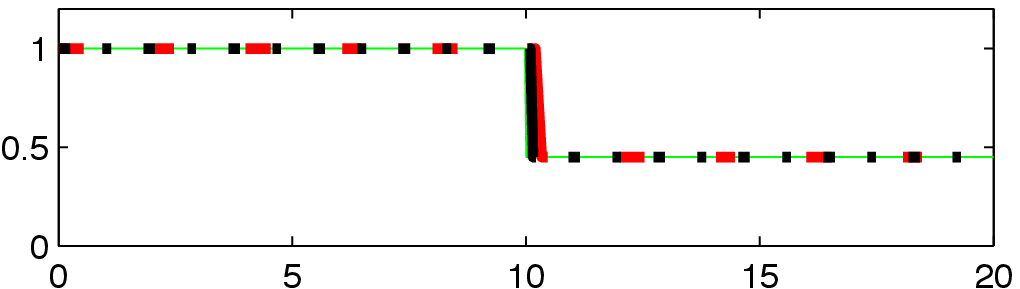}\\
& $x$\\
$t=2.0$\\
\rotatebox{90}{\hspace{10ex}$h$ } &
\includegraphics{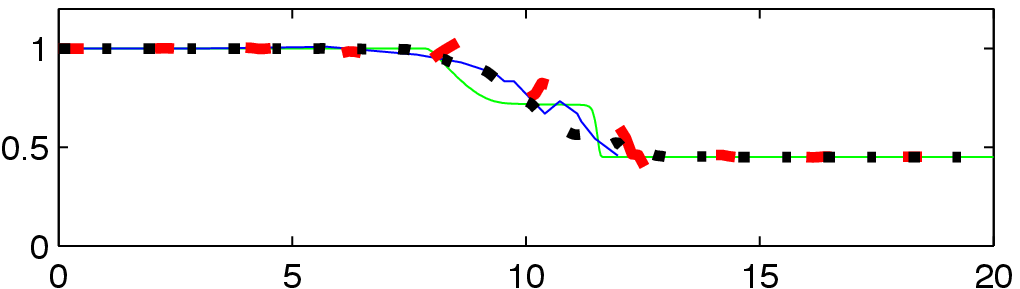}\\
& $x$\\
$t=5.2$\\
\rotatebox{90}{\hspace{10ex}$h$  } &
\includegraphics{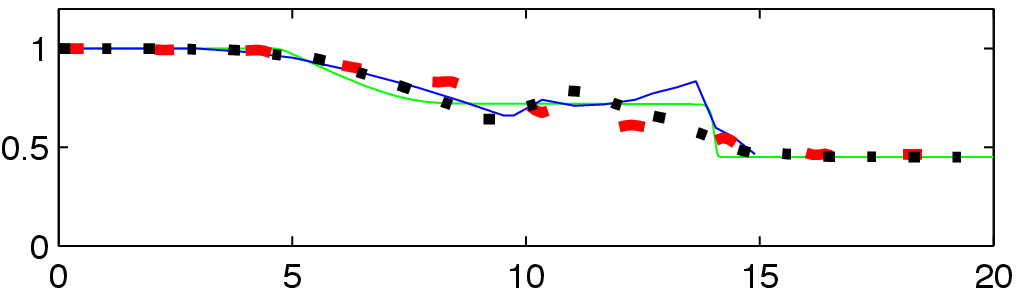}\\
& $x$\\
$t=7.6$\\
\rotatebox{90}{\hspace{10ex}$h$  } &
\includegraphics{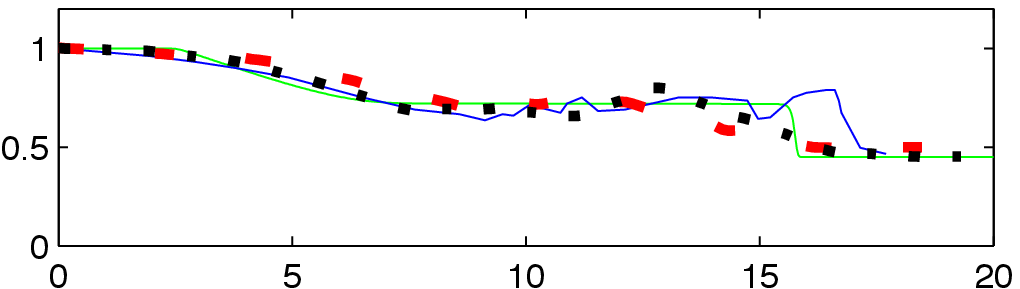}\\
& $x$
\end{tabular}
\caption{Comparison among the simulations of dam breaking: (red)~low resolution gap-tooth with $m=10$ patches; (black)~medium resolution gap-tooth with $m=22$ patches; (green)~the microscale simulation over the whole domain; and (blue) some experimental data \cite[Fig.~8c]{Stansby1998}.
The dam is inside a patch.
The nondimensional shallow depth is $h=0.45$ in front of the dam.
The scale ratio $r=1/6$, and the microscale step $d=1/15$\,.
}
\label{patchdamID}
\end{figure}

\paragraph{The case of the dam being within a patch}
Figure~\ref{patchdamID} plots at four times the gap-tooth simulation (red for $m=10$~patches and black for $m=22$~patches), the microscale simulation over the whole domain~(green), and the experimental data~(blue).
%To avoid singular points at the dam, represent the dam by the function~$\tanh(x)$.
The time $t=0$ graph shows that the water has depth one upstream and depth~$0.45$ downstream, which corresponds to the depth ratio of~$0.45$ in the experiments~\cite[Fig.8c]{Stansby1998}.
Recall the nondimensional length scale~$H$, velocity~$\sqrt{gH}$, and time~$\sqrt{H/g}$.
Nondimensionalise the time in the experiments by $\sqrt{0.1\m/(10\m/\text{s}^{-2})}=0.1$\,s. Therefore the plots of \cite{Stansby1998} are at nondimensional times $t=0,2,5.2$ and~$7.6$. 

The $t=2,5.2,7.6$ graphs in Figure~\ref{patchdamID} show that a turbulent bore forms in all simulations.
However, the turbulent bore in the gap-tooth simulation is smoothed by the relatively large spacing between patches: the medium resolution simulation being noticeably more better than the low resolution.
%At later times $t=5.2$ and $t=7.6$ graphs, the gap-tooth simulation agrees better with the experimental data for the waves behind the water bore.
%The height of the bore in the gap-tooth simulation is smaller than both in the experimental and in the microscale simulation over the whole domain.
The bore in the gap-tooth simulation lags that in the experiment, while the bore in the microscale simulation over the whole domain reasonably tracks that in the experiment until the last time.
The gap-tooth simulation has error $\mathcal O(D^2)\sim0.8$ for the medium resolution macroscale step $D=L/m=20/22\approx0.9$, while the microscale simulation over the whole macroscale domain involves the error~$\mathcal O(d^2)\sim0.001$ for the microscale step $d=2rD/(n+1)=0.03$.
When the number of patches increases, the error~$\mathcal O(D^2)$ decreases, then the gap-tooth simulation performs better, as shown by the black data in comparison to the coarser red data.

However, the gap-tooth simulation saves computer time.
The gap-tooth scheme only takes a compute time of~$0.78$\,s for $m=22$ patches to simulate to the $t=7.6$ graph in Figure~\ref{patchdamID} (all simulations used Matlab with \verb|ode15s| for time integration).
Whereas the microscale simulation over the whole domain with the same microscale step $d=0.03$ needed a compute time of~$74.3$\,s to simulate over the same time.
That is, the whole domain simulation is nearly a hundred times slower than the gap-tooth simulation.
Such a speed-up in this simple pilot study suggests, especially with smaller ratio~\(r\) and implementing projective integration~\cite{Gear02b, Kevrekidis09a}, that the gap-tooth scheme may empower simulation and analysis of large scale problems that are otherwise inaccessible.

\begin{figure}
\centering
\begin{tabular}{c@{\ }c}
\rotatebox{90}{\hspace{15ex} water area} &
\includegraphics[height=45ex]{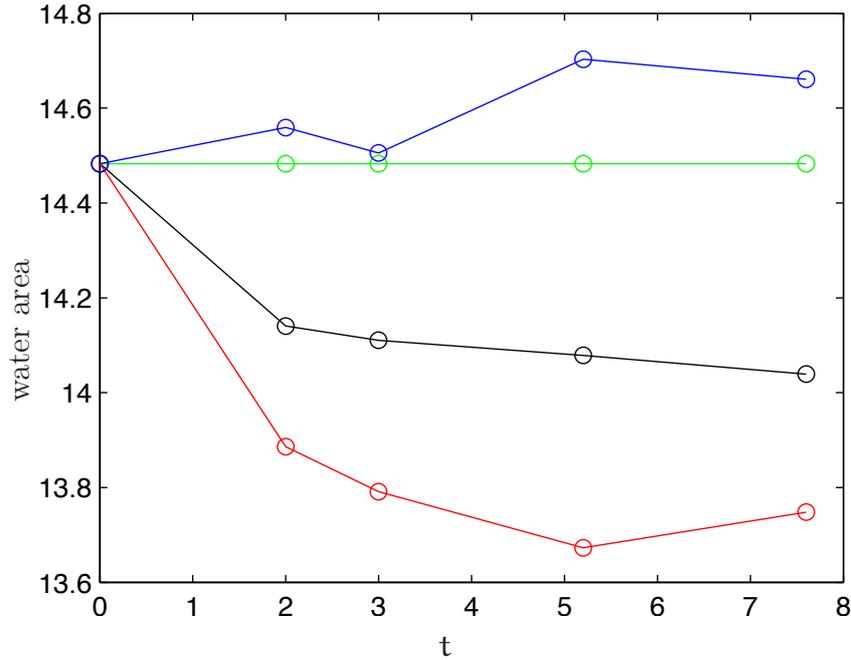}\\
& $t$
\end{tabular}
\caption{The water area over time of the simulations in Figure~\ref{patchdamID}: (black)~gap-tooth with $m=22$ patches; gap-tooth with (red) $m=10$ patches; (green) the microscale simulations over the whole space domain; and (blue)~experiments by \cite{Stansby1998}.
}
\label{patchdamArea}
\end{figure}

Figure~\ref{patchdamArea} shows  in time the water area of the graphs in Figure~\ref{patchdamID}.
Since water is conserved, ideally these curves should be horizontal as seen for the green curve of the microscale simulation over the whole domain.
The red curve shows that the low resolution ($m=10$ patches) gap-tooth simulation loses about~6\% area, mainly in the initial dam break.  The medium resolution simulation loses~2\% (black curve), and again mostly in the initial dam break. 
Both of these loses are due to the relatively coarse spacing of the patches.
The blue curve shows that the experiments gained fluid area, possibly due to entrainment of air in the turbulent bore~\cite[Fig.~8c]{Stansby1998}.
%Initially, they have the same water area.
%The green curve remains unchanged and indicates that the microscale simulation over the whole domain conserves water.
%Stansby's experimental~(blue curve) shows a bigger water area possibly due to the entrainment of air, whereas the gap-tooth simulation~(red and black curves) leads to smaller water area due to the relatively coarse separation of the patches.
%However, when the number of patches~$m$ increase, the water area in the gap-tooth simulation would approach that in the microscale simulation over the whole domain~(green curve).
 
\begin{figure}
\centering
\begin{tabular}{r@{\ }c}
$t=0.0$\\
\rotatebox{90}{\hspace{10ex}$h$} &
\includegraphics{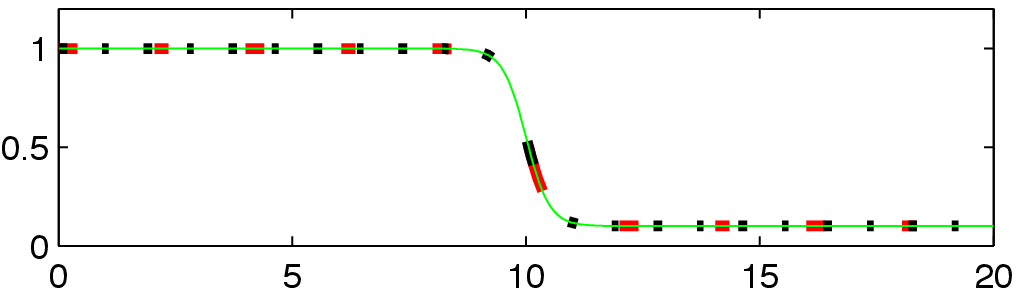}\\
& $x$\\
$t=2.4$\\
\rotatebox{90}{\hspace{10ex}$h$ } &
\includegraphics{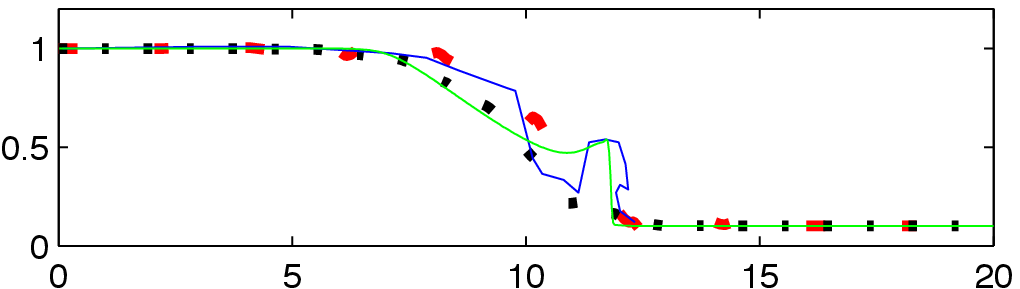}\\
& $x$\\
$t=4.0$\\
\rotatebox{90}{\hspace{10ex}$h$  } &
\includegraphics{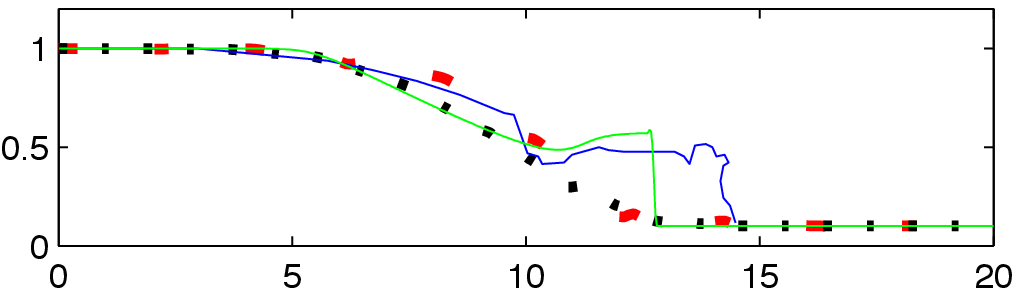}\\
& $x$\\
$t=6.6$\\
\rotatebox{90}{\hspace{10ex}$h$  } &
\includegraphics{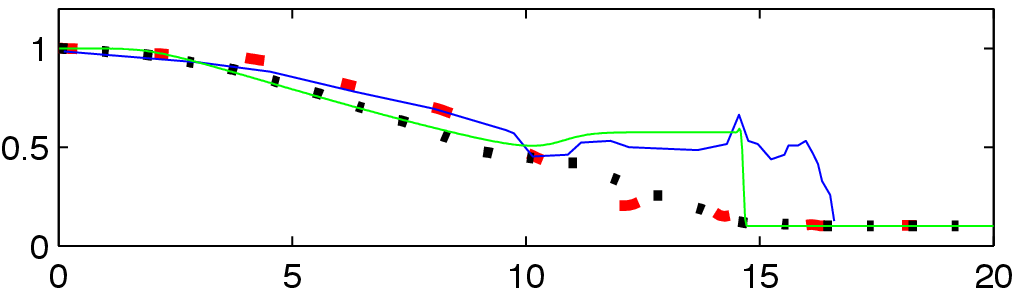}\\
& $x$
\end{tabular}
\caption{Comparison among the simulations of dam-breaking into shallow water of depth $h=0.1$: (red)~low resolution gap-tooth with $m=10$ patches; (black)~medium resolution gap-tooth with $m=22$ patches; (green)~the microscale simulation over the whole domain; and (blue) some experimental data \cite[Fig.~8c]{Stansby1998}.
The patch ratio $r=1/8$, and the microscale step $d=1/20$. 
}
\label{patchdamIDL}
\end{figure}

A further experiment had significantly shallower water in front of the dam which we also simulated.
To avoid singularities being generated in the simulations near the dam, the dam was smoothed to a hyperbolic tangent and the smaller patch ratio \(r=1/8\) was used.
Figure~\ref{patchdamIDL} plots the gap-tooth simulation~(red and black curves for $m=10$ and $m=22$ patches respectively), the microscale simulation over the whole macroscale domain~(green curve), and the experimental data~(blue curve) \cite[Fig.~8b]{Stansby1998}.
The initial shallow depth in front of the dam is~$0.1$ (nondimensional).
Compared with Figure~\ref{patchdamID}, the heights of the turbulent bore in the gap-tooth simulations are significantly smaller than that in the experiment and microscale simulation over the whole domain.
Again, this seems due to the relatively coarse macroscale resolution in these gap-tooth simulations. 
%Because, near the bore, the involved error ratio~$\Ord{D^2}$/\Ord{d^2} between the gap-tooth simulation and microscale simulation over the whole domain becomes big.
The gap-tooth simulation should better approximate the dam-breaking waves with more patches or deeper water in front of the dam.

\begin{figure}
\centering
\begin{tabular}{r@{\ }c}
$t=0.0$\\
\rotatebox{90}{\hspace{10ex}$h$} &
\includegraphics{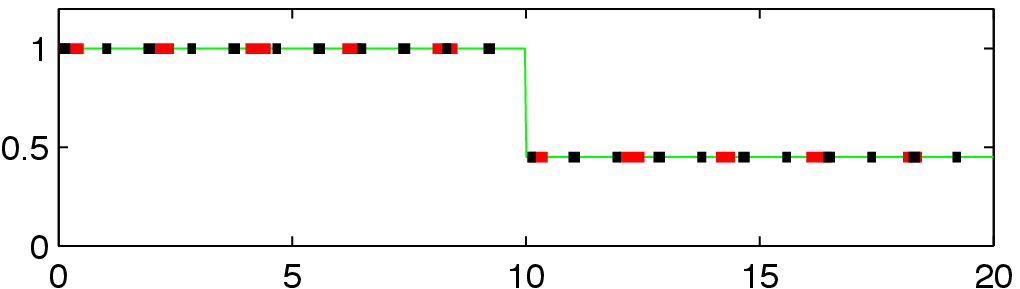}\\
& $x$\\
$t=2.0$\\
\rotatebox{90}{\hspace{10ex}$h$ } &
\includegraphics{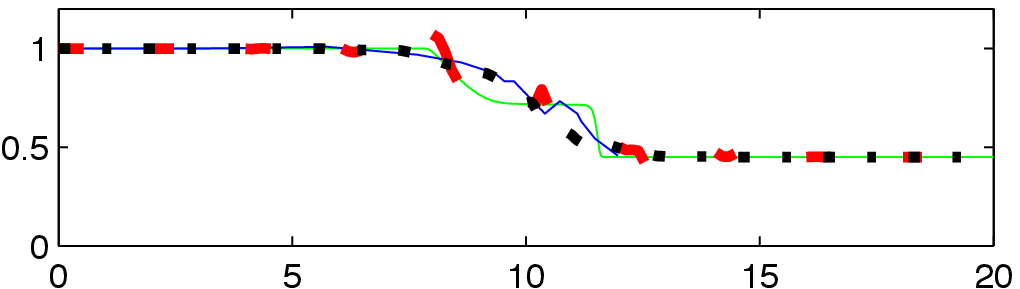}\\
& $x$\\
$t=5.2$\\
\rotatebox{90}{\hspace{10ex}$h$  } &
\includegraphics{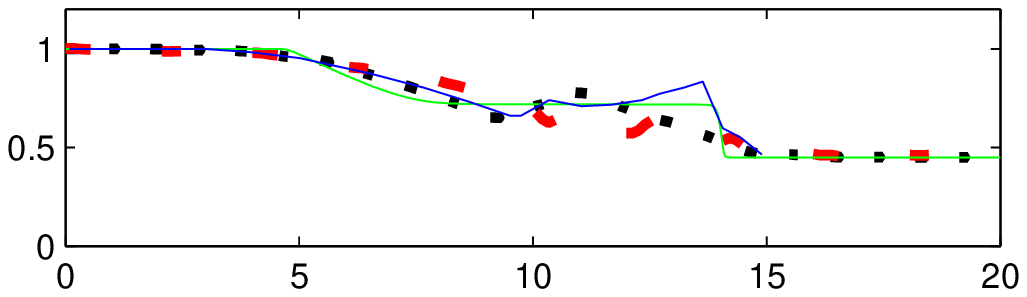}\\
& $x$\\
$t=7.6$\\
\rotatebox{90}{\hspace{10ex}$h$  } &
\includegraphics{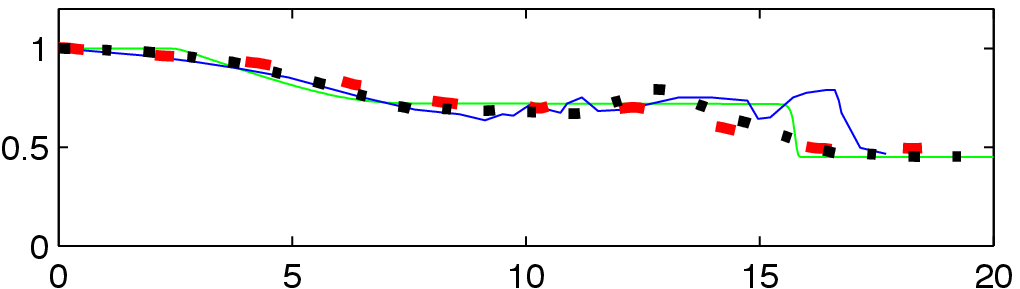}\\
& $x$
\end{tabular}
\caption{Comparison among the simulations of dam breaking: (red)~low resolution gap-tooth with $m=10$ patches; (black)~medium resolution gap-tooth with $m=22$ patches; (green)~the microscale simulation over the whole domain; and (blue)~the experimental data \cite[Fig.~8c]{Stansby1998}.
Initially the dam lies in between patches.
The nondimensional shallow depth is $h=0.45$ in the front of the dam; the patch ratio $r=1/6$; and the microscale step $d=1/15$. 
}
\label{patchdamWD}
\end{figure}

\paragraph{Place the dam between two patches}
Figure~\ref{patchdamWD} plots the gap-tooth simulation~(red and black curves for $m=10$ and $m=22$ patches respectively), the microscale simulation over the whole macroscale domain~(green curve), and the experimental data~(blue curve) by~\cite{Stansby1998} at four times.
The $t=0$ graph shows the initial depth, corresponding to the depth ratio of~$0.45$ in the experiments.
The dam is not resolved within a patch in these gap-tooth simulation.
The $t=2$ graph indicates that this gap-tooth simulation does not appear to be as good as the corresponding results in Figure~\ref{patchdamID}.
Then the $t=5.2$ graph shows there seems to be significantly more microscale oscillations in the low resolution case when compared to Figure~\ref{patchdamID}.
It appears that putting the dam within a patch is better.

\subsection{Boundary conditions for gap-tooth simulation}
\label{patch:damB}

This subsection discusses the invoked boundary conditions at the upstream and downstream boundaries of the macroscale domain in the gap-tooth simulation of the dam-breaking.
Such boundary conditions will be needed in general simulations.

The gap-tooth simulation of the dam-breaking requires boundary conditions at the upstream $x=0$ and downstream $x=L$.
Typically, no-flow boundary conditions are usually implemented in dam breaking \cite[e.g.]{Abdolmaleki2004, Ozgokmen2007a}.
We consider three types of boundary condition at the upstream $x=0$ and downstream $x=L$:
\begin{itemize}
\item constant depth such as $h=1$ at the upstream or depth $h=0.45$ at the downstream;
\item no flux, $\D x\uu=0$, at the upstream or at the downstream;
\item and zero turbulent mean velocity, $\uu=0$, at the upstream or at the downstream, which is equivalent to no fluid flowing through the upstream or downstream, $\D xh=0$ according to the momentum \pde~\eqref{patch:Nu}.
\end{itemize}
Since these boundary conditions could be applied either at the centre or edge of a patch, there could be at least $\binom32^2=36$ possible combinations of boundary conditions in the gap-tooth simulation.

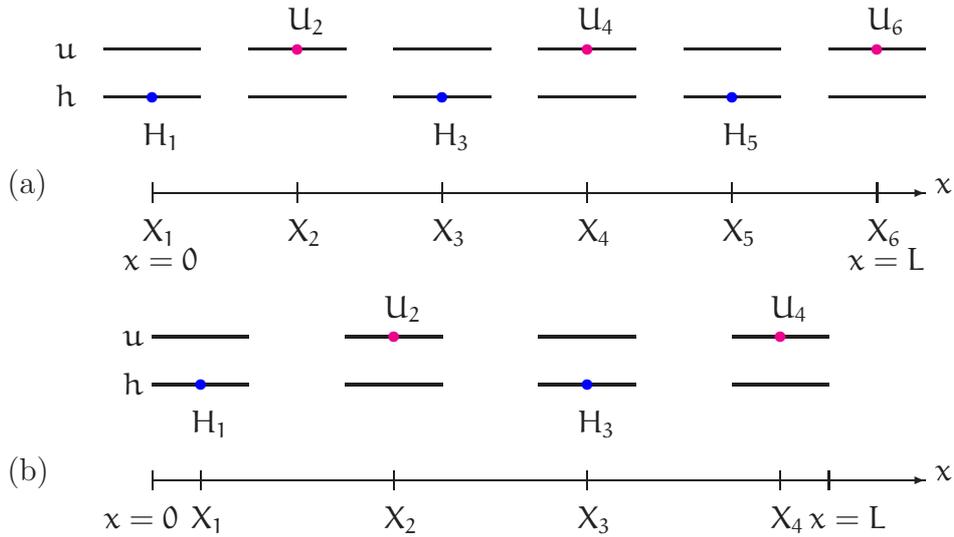
\begin{figure}
\centering
\setlength{\unitlength}{0.7ex}
\begin{picture}(70,60)
{
\put(0,40){\vector(1,0){80}}
\multiput(0,39)(15,0){6}{\line(0,1){2}}
\multiput(-5,50)(15,0){6}{\thicklines\line(1,0){10}}
\multiput(-5,55)(15,0){6}{\thicklines\line(1,0){10}}
\multiput(0,50)(30,0){3}{\color{blue}\circle*{1}}
\multiput(15,55)(30,0){3}{\color{magenta}\circle*{1}}
%\newcounter{i}
\setcounter{i}{0}
\multiput(0,37)(15,0){6}{
  \stepcounter{i}%
  \put(-1,-2){$X_{\arabic{i}}$}}
\put(-1,45){$H_1$} \put(29,45){$H_3$} \put(59,45){$H_5$}
 \put(14,57){$U_2$} \put(44,57){$U_4$} \put(74,57){$U_6$}
\put(-3,32){$x=0$} \put(72,32){$x=L$} \put(81,40){$x$}
\put(-10,49){$h$} \put(-10,54){$u$}
}
\put(-15,40){(a)}
{
\put(0,10){\vector(1,0){80}}
\multiput(5,9)(20,0){4}{\line(0,1){2}}
\multiput(0,20)(20,0){4}{\thicklines\line(1,0){10}}
\multiput(0,25)(20,0){4}{\thicklines\line(1,0){10}}
\multiput(5,20)(40,0){2}{\color{blue}\circle*{1}}
\multiput(25,25)(40,0){2}{\color{magenta}\circle*{1}}
\setcounter{i}{0}
\multiput(5,7)(20,0){4}{
  \stepcounter{i}%
  \put(-1,-2){$X_{\arabic{i}}$}}
\put(4,15){$H_1$} \put(44,15){$H_3$} \put(24,27){$U_2$} \put(64,27){$U_4$}  
 \put(81,10){$x$} \put(0,9){\line(0,1){2}\put(-5,-4){$x=0$}} \put(70,9){\line(0,1){2}\put(-2,-4){$x=L$}}
 \put(-3,19){$h$} \put(-3,24){$u$}
}
\put(-15,10){(b)}
\end{picture}
\caption{There are at least two ways to implement boundary conditions at the upstream $x=0$ and downstream $x=L$ in gap-tooth simulation:
(a)/(b),~respectively, invokes the boundary conditions for the macroscale\slash microscale values of the leftmost and rightmost patches.}
\label{patchDambc}
\end{figure}%

In the gap-tooth simulation, boundary conditions are invoked to either the macroscale or microscale values on the leftmost and rightmost patches, as shown schematically by Figure~\ref{patchDambc}. 
Figure~\ref{patchdamID}--\ref{patchdamWD} implement the boundary conditions $\uu_{1,1}=0$ as drawn in Figure~\ref{patchDambc}(b), and $h_{m,1}=h_{m,n+2}=H_{m-1}$, where $\uu_{1,1}$~is the left edge of the first patch, $h_{m,1}$ and~$h_{m,n+2}$ are the edges of the $m$th~patch, and $H_{m-1}$~is the macroscale value on the $(m-1)$th~patch through the coupling conditions~\eqref{patch:cph} with the assumption of zero values on the fictitious $(m+1)$th patch in the simulation.
Further work could explore the gap-tooth simulation with different boundary conditions.

\section{Conclusion}
\label{patch:conclusion}

Developing some preliminary research \cite[]{Cao2013}, we explored the gap-tooth scheme both theoretically and with a highly nonlinear microscale simulator~\eqref{eqs:patch:Nji} of turbulent shallow water waves (Section~\ref{sec:micro}).
The resultant numerical simulations indicated that the gap-tooth scheme on a staggered macroscale grid is useful for wave-like systems.
Section~\ref{sec:Eig} reported numerical eigenvalue analysis that clearly showed, in Figure~\ref{patchNeig}, the separation between relatively slow macroscale wave modes and the microscale fast waves supported within patches. 
The theoretical support of section~\ref{sec:consistency} proves for a wide range of dispersive linear wave-like systems that the gap-tooth scheme generates macroscale simulations consistent with the microscale. 
Such consistency holds for a much wider class of nonlinear wave-like systems (section~\ref{sec:consistencyC}).
Section~\ref{sec:nsm} establishes that the gap-tooth scheme has a sound closure in terms of macroscale variables, but with the caveat that resonance among significant microscale waves could cause differing macroscale simulations.
Section~\ref{sec:dam} applied the gap-tooth scheme to the highly nonlinear flow of dam-breaking waves.
Figure~\ref{patchdamID}--\ref{patchdamWD} shows that although the turbulent bore lags and the height of this bore is a bit smaller, we reasonably predict the dam-breaking.
THe major limitation in the gap-tooth scheme appears to be that it cannot resolve microscale dynamics between the patches so it is primarily useful for macroscale dynamics which are globally varying on the macroscale.
Nonetheless, there is scope for resolving microscale dynamics by, for example, putting the rapid changes associated with the dam break within a patch (section~\ref{patch:damN}).

\paragraph{Acknowledgements}
Part of this research was supported by grant DP120104260 from the Australian Research Council.
We thank Prof.~Yannis Kevrekidis for inspiring discussions.

\bibliographystyle{agsm}
\bibliography{bibexport}

\appendix
\section{Ancillary computer algebra program}
\label{sec:appendix}

This computer algebra code constructs the slow manifolds for the gap-tooth scheme discussed in Section~\ref{sec:consistencyC}.
We use the Reduce computer algebra package\footnote{\url{http://www.reduce-algebra.com/}} because it is freely available and because it is perhaps the fastest general purpose computer algebra system \cite[e.g.]{Fateman2002}.

\small
\verbatimlisting{caConsistencyH.txt}

\end{document}